\ifx\shlhetal\undefinedcontrolsequence\let\shlhetal\relax\fi

\input amstex
\expandafter\ifx\csname mathdefs.tex\endcsname\relax
  \expandafter\gdef\csname mathdefs.tex\endcsname{}
\else \message{Hey!  Apparently you were trying to
  \string\input{mathdefs.tex} twice.   This does not make sense.} 
\errmessage{Please edit your file (probably \jobname.tex) and remove
any duplicate ``\string\input'' lines}\endinput\fi




\catcode`\X=12\catcode`\@=11

\def\n@wcount{\alloc@0\count\countdef\insc@unt}
\def\n@wwrite{\alloc@7\write\chardef\sixt@@n}
\def\n@wread{\alloc@6\read\chardef\sixt@@n}
\def\r@s@t{\relax}\def\v@idline{\par}\def\@mputate#1/{#1}
\def\l@c@l#1X{\firstpart.#1}\def\gl@b@l#1X{#1}\def\t@d@l#1X{{}}

\def\crossrefs#1{\ifx\all#1\let\tr@ce=\all\else\def\tr@ce{#1,}\fi
   \n@wwrite\cit@tionsout\openout\cit@tionsout=\jobname.cit 
   \write\cit@tionsout{\tr@ce}\expandafter\setfl@gs\tr@ce,}
\def\setfl@gs#1,{\def\@{#1}\ifx\@\empty\let\next=\relax
   \else\let\next=\setfl@gs\expandafter\xdef
   \csname#1tr@cetrue\endcsname{}\fi\next}
\def\m@ketag#1#2{\expandafter\n@wcount\csname#2tagno\endcsname
     \csname#2tagno\endcsname=0\let\tail=\all\xdef\all{\tail#2,}
   \ifx#1\l@c@l\let\tail=\r@s@t\xdef\r@s@t{\csname#2tagno\endcsname=0\tail}\fi
   \expandafter\gdef\csname#2cite\endcsname##1{\expandafter
     \ifx\csname#2tag##1\endcsname\relax?\else\csname#2tag##1\endcsname\fi
     \expandafter\ifx\csname#2tr@cetrue\endcsname\relax\else
     \write\cit@tionsout{#2tag ##1 cited on page \folio.}\fi}
   \expandafter\gdef\csname#2page\endcsname##1{\expandafter
     \ifx\csname#2page##1\endcsname\relax?\else\csname#2page##1\endcsname\fi
     \expandafter\ifx\csname#2tr@cetrue\endcsname\relax\else
     \write\cit@tionsout{#2tag ##1 cited on page \folio.}\fi}
   \expandafter\gdef\csname#2tag\endcsname##1{\expandafter
      \ifx\csname#2check##1\endcsname\relax
      \expandafter\xdef\csname#2check##1\endcsname{}%
      \else\immediate\write16{Warning: #2tag ##1 used more than once.}\fi
      \multit@g{#1}{#2}##1/X%
      \write\t@gsout{#2tag ##1 assigned number \csname#2tag##1\endcsname\space
      on page \number\count0.}%
   \csname#2tag##1\endcsname}}

\def\multit@g#1#2#3/#4X{\def\t@mp{#4}\ifx\t@mp\empty%
      \global\advance\csname#2tagno\endcsname by 1 
      \expandafter\xdef\csname#2tag#3\endcsname
      {#1\number\csname#2tagno\endcsnameX}%
   \else\expandafter\ifx\csname#2last#3\endcsname\relax
      \expandafter\n@wcount\csname#2last#3\endcsname
      \global\advance\csname#2tagno\endcsname by 1 
      \expandafter\xdef\csname#2tag#3\endcsname
      {#1\number\csname#2tagno\endcsnameX}
      \write\t@gsout{#2tag #3 assigned number \csname#2tag#3\endcsname\space
      on page \number\count0.}\fi
   \global\advance\csname#2last#3\endcsname by 1
   \def\t@mp{\expandafter\xdef\csname#2tag#3/}%
   \expandafter\t@mp\@mputate#4\endcsname
   {\csname#2tag#3\endcsname\lastpart{\csname#2last#3\endcsname}}\fi}
\def\t@gs#1{\def\all{}\m@ketag#1e\m@ketag#1s\m@ketag\t@d@l p
\let\realscite\scite
\let\realstag\stag
   \m@ketag\gl@b@l r \n@wread\t@gsin
   \openin\t@gsin=\jobname.tgs \re@der \closein\t@gsin
   \n@wwrite\t@gsout\openout\t@gsout=\jobname.tgs }
\outer\def\localtags{\t@gs\l@c@l}
\outer\def\globaltags{\t@gs\gl@b@l}
\outer\def\newlocaltag#1{\m@ketag\l@c@l{#1}}
\outer\def\newglobaltag#1{\m@ketag\gl@b@l{#1}}

\newif\ifpr@ 
\def\m@kecs #1tag #2 assigned number #3 on page #4.%
   {\expandafter\gdef\csname#1tag#2\endcsname{#3}
   \expandafter\gdef\csname#1page#2\endcsname{#4}
   \ifpr@\expandafter\xdef\csname#1check#2\endcsname{}\fi}
\def\re@der{\ifeof\t@gsin\let\next=\relax\else
   \read\t@gsin to\t@gline\ifx\t@gline\v@idline\else
   \expandafter\m@kecs \t@gline\fi\let \next=\re@der\fi\next}
\def\pretags#1{\pr@true\pret@gs#1,,}
\def\pret@gs#1,{\def\@{#1}\ifx\@\empty\let\n@xtfile=\relax
   \else\let\n@xtfile=\pret@gs \openin\t@gsin=#1.tgs \message{#1} \re@der 
   \closein\t@gsin\fi \n@xtfile}

\newcount\sectno\sectno=0\newcount\subsectno\subsectno=0
\newif\ifultr@local \def\ultralocal{\ultr@localtrue}
\def\firstpart{\number\sectno}
\def\lastpart#1{\ifcase#1 \or a\or b\or c\or d\or e\or f\or g\or h\or 
   i\or k\or l\or m\or n\or o\or p\or q\or r\or s\or t\or u\or v\or w\or 
   x\or y\or z \fi}

\def\resetall{\global\advance\sectno by 1\subsectno=0
   \gdef\firstpart{\number\sectno}\r@s@t}
\def\resetsub{\global\advance\subsectno by 1
   \gdef\firstpart{\number\sectno.\number\subsectno}\r@s@t}
\def\newsection#1\par{\resetall\vskip0pt plus.3\vsize\penalty-250
   \vskip0pt plus-.3\vsize\bigskip\bigskip
   \message{#1}\leftline{\bf#1}\nobreak\bigskip}
\def\subsection#1\par{\ifultr@local\resetsub\fi
   \vskip0pt plus.2\vsize\penalty-250\vskip0pt plus-.2\vsize
   \bigskip\smallskip\message{#1}\leftline{\bf#1}\nobreak\medskip}


\newdimen\marginshift

\newdimen\margindelta
\newdimen\marginmax
\newdimen\marginmin

\def\margininit{       
\marginmax=3 true cm                  
				      
\margindelta=0.1 true cm              
\marginmin=0.1true cm                 
\marginshift=\marginmin
}    

\def\t@gsjj#1,{\def\@{#1}\ifx\@\empty\let\next=\relax\else\let\next=\t@gsjj
   \def\@@{p}\ifx\@\@@\else
   \expandafter\gdef\csname#1cite\endcsname##1{\citejj{##1}}
   \expandafter\gdef\csname#1page\endcsname##1{?}
   \expandafter\gdef\csname#1tag\endcsname##1{\tagjj{##1}}\fi\fi\next}
\newif\ifshowstuffinmargin
\showstuffinmarginfalse
\def\jjtags{\ifx\shlhetal\relax 
  \else
\ifx\shlhetal\undefinedcontrolseq
\else
\showstuffinmargintrue
\ifx\all\relax\else\expandafter\t@gsjj\all,\fi\fi \fi
}

\def\tagjj#1{\realstag{#1}\oldmginpar{\zeigen{#1}}}
\def\citejj#1{\rechnen{#1}\mginpar{\zeigen{#1}}}     

\def\rechnen#1{\expandafter\ifx\csname stag#1\endcsname\relax ??\else
                           \csname stag#1\endcsname\fi}

\newdimen\theight

\def\marginfont{\sevenrm}

\def\trymarginbox#1{\setbox0=\hbox{\marginfont\hskip\marginshift #1}%
		\global\marginshift\wd0 
		\global\advance\marginshift\margindelta}

\def \oldmginpar#1{%
\ifvmode\setbox0\hbox to \hsize{\hfill\rlap{\marginfont\quad#1}}%
\ht0 0cm
\dp0 0cm
\box0\vskip-\baselineskip
\else 
             \vadjust{\trymarginbox{#1}%
		\ifdim\marginshift>\marginmax \global\marginshift\marginmin
			\trymarginbox{#1}%
                \fi
             \theight=\ht0
             \advance\theight by \dp0    \advance\theight by \lineskip
             \kern -\theight \vbox to \theight{\rightline{\rlap{\box0}}%
\vss}}\fi}

\newdimen\upordown
\global\upordown=8pt
\font\tinyfont=cmtt8 
\def\mginpar#1{\smash{\hbox to 0cm{\kern-10pt\raise7pt\hbox{\tinyfont #1}\hss}}}
\def\mginpar#1{{\hbox to 0cm{\kern-10pt\raise\upordown\hbox{\tinyfont #1}\hss}}\global\upordown-\upordown}


\def\t@gsoff#1,{\def\@{#1}\ifx\@\empty\let\next=\relax\else\let\next=\t@gsoff
   \def\@@{p}\ifx\@\@@\else
   \expandafter\gdef\csname#1cite\endcsname##1{\zeigen{##1}}
   \expandafter\gdef\csname#1page\endcsname##1{?}
   \expandafter\gdef\csname#1tag\endcsname##1{\zeigen{##1}}\fi\fi\next}
\def\verbatimtags{\showstuffinmarginfalse
\ifx\all\relax\else\expandafter\t@gsoff\all,\fi}
\def\zeigen#1{\hbox{$\scriptstyle\langle$}#1\hbox{$\scriptstyle\rangle$}}


\def\margintag#1{\ifshowstuffinmargin\oldmginpar{\zeigen{#1}}\fi}

\def\(#1){\edef\dot@g{\ifmmode\ifinner(\hbox{\noexpand\etag{#1}})
   \else\noexpand\eqno(\hbox{\noexpand\etag{#1}})\fi
   \else(\noexpand\ecite{#1})\fi}\dot@g}

\newif\ifbr@ck
\def\eat#1{}
\def\[#1]{\br@cktrue[\br@cket#1'X]}
\def\br@cket#1'#2X{\def\temp{#2}\ifx\temp\empty\let\next\eat
   \else\let\next\br@cket\fi
   \ifbr@ck\br@ckfalse\br@ck@t#1,X\else\br@cktrue#1\fi\next#2X}
\def\br@ck@t#1,#2X{\def\temp{#2}\ifx\temp\empty\let\neext\eat
   \else\let\neext\br@ck@t\def\temp{,}\fi
   \def\teemp{#1}\ifx\teemp\empty\else\rcite{#1}\fi\temp\neext#2X}
\def\resetbr@cket{\gdef\[##1]{[\rtag{##1}]}}
\def\references{\resetbr@cket\newsection References\par}

\newtoks\symb@ls\newtoks\s@mb@ls\newtoks\p@gelist\n@wcount\ftn@mber
    \ftn@mber=1\newif\ifftn@mbers\ftn@mbersfalse\newif\ifbyp@ge\byp@gefalse
\def\defm@rk{\ifftn@mbers\n@mberm@rk\else\symb@lm@rk\fi}
\def\n@mberm@rk{\xdef\m@rk{{\the\ftn@mber}}%
    \global\advance\ftn@mber by 1 }
\def\rot@te#1{\let\temp=#1\global#1=\expandafter\r@t@te\the\temp,X}
\def\r@t@te#1,#2X{{#2#1}\xdef\m@rk{{#1}}}
\def\b@@st#1{{$^{#1}$}}\def\str@p#1{#1}
\def\symb@lm@rk{\ifbyp@ge\rot@te\p@gelist\ifnum\expandafter\str@p\m@rk=1 
    \s@mb@ls=\symb@ls\fi\write\f@nsout{\number\count0}\fi \rot@te\s@mb@ls}
\def\byp@ge{\byp@getrue\n@wwrite\f@nsin\openin\f@nsin=\jobname.fns 
    \n@wcount\currentp@ge\currentp@ge=0\p@gelist={0}
    \re@dfns\closein\f@nsin\rot@te\p@gelist
    \n@wread\f@nsout\openout\f@nsout=\jobname.fns }
\def\m@kelist#1X#2{{#1,#2}}
\def\re@dfns{\ifeof\f@nsin\let\next=\relax\else\read\f@nsin to \f@nline
    \ifx\f@nline\v@idline\else\let\t@mplist=\p@gelist
    \ifnum\currentp@ge=\f@nline
    \global\p@gelist=\expandafter\m@kelist\the\t@mplistX0
    \else\currentp@ge=\f@nline
    \global\p@gelist=\expandafter\m@kelist\the\t@mplistX1\fi\fi
    \let\next=\re@dfns\fi\next}
\def\symbols#1{\symb@ls={#1}\s@mb@ls=\symb@ls} 
\def\bigsymbol{\textstyle}
\symbols{\bigsymbol\ast,\dagger,\ddagger,\sharp,\flat,\natural,\star}
\def\ftnumbers{\ftn@mberstrue} \def\ftsymbols{\ftn@mbersfalse}
\def\paginal{\byp@ge} \def\resetftnumbers{\ftn@mber=1}
\def\ftnote#1{\defm@rk\expandafter\expandafter\expandafter\footnote
    \expandafter\b@@st\m@rk{#1}}

\long\def\jump#1\endjump{}
\def\ssum{\mathop{\lower .1em\hbox{$\textstyle\Sigma$}}\nolimits}

\def\qed{\nobreak\kern 1em \vrule height .5em width .5em depth 0em}
\def\newneq{\hbox{\rlap{\hbox to 1\wd9{\hss$=$\hss}}\raise .1em 
   \hbox to 1\wd9{\hss$\scriptscriptstyle/$\hss}}}
\def\subsetne{\setbox9 = \hbox{$\subset$}\mathrel{\hbox{\rlap
   {\lower .4em \newneq}\raise .13em \hbox{$\subset$}}}}
\def\supsetne{\setbox9 = \hbox{$\subset$}\mathrel{\hbox{\rlap
   {\lower .4em \newneq}\raise .13em \hbox{$\supset$}}}}

\def\vbar{\mathchoice{\vrule height6.3ptdepth-.5ptwidth.8pt\kern-.8pt}
   {\vrule height6.3ptdepth-.5ptwidth.8pt\kern-.8pt}
   {\vrule height4.1ptdepth-.35ptwidth.6pt\kern-.6pt}
   {\vrule height3.1ptdepth-.25ptwidth.5pt\kern-.5pt}}
\def\f@dge{\mathchoice{}{}{\mkern.5mu}{\mkern.8mu}}
\def\b@c#1#2{{\rm \mkern#2mu\vbar\mkern-#2mu#1}}
\def\b@b#1{{\rm I\mkern-3.5mu #1}}
\def\b@a#1#2{{\rm #1\mkern-#2mu\f@dge #1}}
\def\bb#1{{\count4=`#1 \advance\count4by-64 \ifcase\count4\or\b@a A{11.5}\or
   \b@b B\or\b@c C{5}\or\b@b D\or\b@b E\or\b@b F \or\b@c G{5}\or\b@b H\or
   \b@b I\or\b@c J{3}\or\b@b K\or\b@b L \or\b@b M\or\b@b N\or\b@c O{5} \or
   \b@b P\or\b@c Q{5}\or\b@b R\or\b@a S{8}\or\b@a T{10.5}\or\b@c U{5}\or
   \b@a V{12}\or\b@a W{16.5}\or\b@a X{11}\or\b@a Y{11.7}\or\b@a Z{7.5}\fi}}

\catcode`\X=11 \catcode`\@=12




\let\thischap\jobname

\def\partof#1{\csname returnthe#1part\endcsname}
\def\chapof#1{\csname returnthe#1chap\endcsname}

\def\setchapter#1,#2,#3;{%
  \expandafter\def\csname returnthe#1part\endcsname{#2}%
  \expandafter\def\csname returnthe#1chap\endcsname{#3}%
}

\setchapter 300a,A,II.A;
\setchapter 300b,A,II.B;
\setchapter 300c,A,II.C;
\setchapter 300d,A,II.D;
\setchapter 300e,A,II.E;
\setchapter 300f,A,II.F;
\setchapter 300g,A,II.G;
\setchapter  E53,B,N;
\setchapter  88r,B,I;
\setchapter  600,B,III;
\setchapter  705,B,IV;
\setchapter  734,B,V;

\def\cprefix#1{
\edef\theotherpart{\partof{#1}}\edef\theotherchap{\chapof{#1}}%
\ifx\theotherpart\thispart
   \ifx\theotherchap\thischap 
    \else 
     \theotherchap%
    \fi
   \else 
     \theotherchap\fi}

\def\sectioncite[#1]#2{%
     \cprefix{#2}#1}

\edef\thispart{\partof{\thischap}}
\edef\thischap{\chapof{\thischap}}

\def\lastpage of '#1' is #2.{\expandafter\def\csname lastpage#1\endcsname{#2}}


\def\spuriousreset{}


\expandafter\ifx\csname citeadd.tex\endcsname\relax
\expandafter\gdef\csname citeadd.tex\endcsname{}
\else \message{Hey!  Apparently you were trying to
\string\input{citeadd.tex} twice.   This does not make sense.} 
\errmessage{Please edit your file (probably \jobname.tex) and remove
any duplicate ``\string\input'' lines}\endinput\fi

\sectno=-1   
\localtags
\jjtags
\NoBlackBoxes
\define\mr{\medskip\roster}
\define\sn{\smallskip\noindent}
\define\mn{\medskip\noindent}
\define\bn{\bigskip\noindent}
\define\ub{\underbar}
\define\wilog{\text{without loss of generality}}
\define\ermn{\endroster\medskip\noindent}
\define\dbca{\dsize\bigcap}
\define\dbcu{\dsize\bigcup}
\define \nl{\newline}

\magnification=\magstep 1
\documentstyle{amsppt}

{    
\catcode`@11

\ifx\alicetwothousandloaded@\relax
  \endinput\else\global\let\alicetwothousandloaded@\relax\fi

\gdef\subjclass{\let\savedef@\subjclass
 \def\subjclass##1\endsubjclass{\let\subjclass\savedef@
   \toks@{\def\usualspace{{\rm\enspace}}\eightpoint}%
   \toks@@{##1\unskip.}%
   \edef\thesubjclass@{\the\toks@
     \frills@{{\noexpand\rm2000 {\noexpand\it Mathematics Subject
       Classification}.\noexpand\enspace}}%
     \the\toks@@}}%
  \nofrillscheck\subjclass}
} 


\expandafter\ifx\csname alice2jlem.tex\endcsname\relax
  \expandafter\xdef\csname alice2jlem.tex\endcsname{\the\catcode`@}
\else \message{Hey!  Apparently you were trying to
\string\input{alice2jlem.tex}  twice.   This does not make sense.}
\errmessage{Please edit your file (probably \jobname.tex) and remove
any duplicate ``\string\input'' lines}\endinput\fi

\expandafter\ifx\csname bib4plain.tex\endcsname\relax
  \expandafter\gdef\csname bib4plain.tex\endcsname{}
\else \message{Hey!  Apparently you were trying to \string\input
  bib4plain.tex twice.   This does not make sense.}
\errmessage{Please edit your file (probably \jobname.tex) and remove
any duplicate ``\string\input'' lines}\endinput\fi

\def\renewcommand{\newcommand}	       
\edef\cite{\the\catcode`@}%
\catcode`@ = 11
\let\@oldatcatcode = \cite
\chardef\@letter = 11
\chardef\@other = 12
%
%
%
%
\def\@innerdef#1#2{\edef#1{\expandafter\noexpand\csname #2\endcsname}}%
%
%
\@innerdef\@innernewcount{newcount}%
\@innerdef\@innernewdimen{newdimen}%
\@innerdef\@innernewif{newif}%
\@innerdef\@innernewwrite{newwrite}%
%
%
%
\def\@gobble#1{}%
%
%
%
\ifx\inputlineno\@undefined
   \let\@linenumber = \empty 
\else
   \def\@linenumber{\the\inputlineno:\space}%
\fi
%
%
%
\def\@futurenonspacelet#1{\def\cs{#1}%
   \afterassignment\@stepone\let\@nexttoken=
}%
\begingroup 
\def\\{\global\let\@stoken= }%
\\ 
\endgroup
\def\@stepone{\expandafter\futurelet\cs\@steptwo}%
\def\@steptwo{\expandafter\ifx\cs\@stoken\let\@@next=\@stepthree
   \else\let\@@next=\@nexttoken\fi \@@next}%
\def\@stepthree{\afterassignment\@stepone\let\@@next= }%
%
%
%
\def\@getoptionalarg#1{%
   \let\@optionaltemp = #1%
   \let\@optionalnext = \relax
   \@futurenonspacelet\@optionalnext\@bracketcheck
}%
%
%
\def\@bracketcheck{%
   \ifx [\@optionalnext
      \expandafter\@@getoptionalarg
   \else
      \let\@optionalarg = \empty
      \expandafter\@optionaltemp
   \fi
}%
\def\@@getoptionalarg[#1]{%
   \def\@optionalarg{#1}%
   \@optionaltemp
}%
%
%
%
\def\@nnil{\@nil}%
\def\@fornoop#1\@@#2#3{}%
\def\@for#1:=#2\do#3{%
   \edef\@fortmp{#2}%
   \ifx\@fortmp\empty \else
      \expandafter\@forloop#2,\@nil,\@nil\@@#1{#3}%
   \fi
}%
\def\@forloop#1,#2,#3\@@#4#5{\def#4{#1}\ifx #4\@nnil \else
       #5\def#4{#2}\ifx #4\@nnil \else#5\@iforloop #3\@@#4{#5}\fi\fi
}%
\def\@iforloop#1,#2\@@#3#4{\def#3{#1}\ifx #3\@nnil
       \let\@nextwhile=\@fornoop \else
      #4\relax\let\@nextwhile=\@iforloop\fi\@nextwhile#2\@@#3{#4}%
}%
%
%
%
\@innernewif\if@fileexists
\def\@testfileexistence{\@getoptionalarg\@finishtestfileexistence}%
\def\@finishtestfileexistence#1{%
   \begingroup
      \def\extension{#1}%
      \immediate\openin0 =
         \ifx\@optionalarg\empty\jobname\else\@optionalarg\fi
         \ifx\extension\empty \else .#1\fi
         \space
      \ifeof 0
         \global\@fileexistsfalse
      \else
         \global\@fileexiststrue
      \fi
      \immediate\closein0
   \endgroup
}%
%
%
%
%
\def\bibliographystyle#1{%
   \@readauxfile
   \@writeaux{\string\bibstyle{#1}}%
}%
\let\bibstyle = \@gobble
%
%
\let\bblfilebasename = \jobname
\def\bibliography#1{%
   \@readauxfile
   \@writeaux{\string\bibdata{#1}}%
   \@testfileexistence[\bblfilebasename]{bbl}%
   \if@fileexists
      \nobreak
      \@readbblfile
   \fi
}%
\let\bibdata = \@gobble
%
%
\def\nocite#1{%
   \@readauxfile
   \@writeaux{\string\citation{#1}}%
}%
\@innernewif\if@notfirstcitation
%
%
\def\cite{\@getoptionalarg\@cite}%
%
%
\def\@cite#1{%
   \let\@citenotetext = \@optionalarg
   \printcitestart
   \nocite{#1}%
   \@notfirstcitationfalse
   \@for \@citation :=#1\do
   {%
      \expandafter\@onecitation\@citation\@@
   }%
   \ifx\empty\@citenotetext\else
      \printcitenote{\@citenotetext}%
   \fi
   \printcitefinish
}%
\newif\ifweareinprivate
\weareinprivatetrue
\ifx\shlhetal\undefinedcontrolseq\weareinprivatefalse\fi
\ifx\shlhetal\relax\weareinprivatefalse\fi
\def\@onecitation#1\@@{%
   \if@notfirstcitation
      \printbetweencitations
   \fi
   \expandafter \ifx \csname\@citelabel{#1}\endcsname \relax
      \if@citewarning
         \message{\@linenumber Undefined citation `#1'.}%
      \fi
     \ifweareinprivate
      \expandafter\gdef\csname\@citelabel{#1}\endcsname{%
\strut 
\vadjust{\vskip-\dp\strutbox
\vbox to 0pt{\vss\parindent0cm \leftskip=\hsize 
\advance\leftskip3mm
\advance\hsize 4cm\strut\openup-4pt 
\rightskip 0cm plus 1cm minus 0.5cm ?  #1 ?\strut}}
         {\tt
            \escapechar = -1
            \nobreak\hskip0pt\pfeilsw
            \expandafter\string\csname#1\endcsname
             \pfeilso
            \nobreak\hskip0pt
         }%
      }%
     \else  
      \expandafter\gdef\csname\@citelabel{#1}\endcsname{%
            {\tt\expandafter\string\csname#1\endcsname}
      }%
     \fi  
   \fi
   \csname\@citelabel{#1}\endcsname
   \@notfirstcitationtrue
}%
%
%
\def\@citelabel#1{b@#1}%
%
%
\def\@citedef#1#2{\expandafter\gdef\csname\@citelabel{#1}\endcsname{#2}}%
%
%
%
\def\@readbblfile{%
   \ifx\@itemnum\@undefined
      \@innernewcount\@itemnum
   \fi
   \begingroup
      \def\begin##1##2{%
         \setbox0 = \hbox{\biblabelcontents{##2}}%
         \biblabelwidth = \wd0
      }%
      \def\end##1{}
      %
      %
      \@itemnum = 0
      \def\bibitem{\@getoptionalarg\@bibitem}%
      \def\@bibitem{%
         \ifx\@optionalarg\empty
            \expandafter\@numberedbibitem
         \else
            \expandafter\@alphabibitem
         \fi
      }%
      \def\@alphabibitem##1{%
         \expandafter \xdef\csname\@citelabel{##1}\endcsname {\@optionalarg}%
         \ifx\biblabelprecontents\@undefined
            \let\biblabelprecontents = \relax
         \fi
         \ifx\biblabelpostcontents\@undefined
            \let\biblabelpostcontents = \hss
         \fi
         \@finishbibitem{##1}%
      }%
      \def\@numberedbibitem##1{%
         \advance\@itemnum by 1
         \expandafter \xdef\csname\@citelabel{##1}\endcsname{\number\@itemnum}%
         \ifx\biblabelprecontents\@undefined
            \let\biblabelprecontents = \hss
         \fi
         \ifx\biblabelpostcontents\@undefined
            \let\biblabelpostcontents = \relax
         \fi
         \@finishbibitem{##1}%
      }%
      \def\@finishbibitem##1{%
         \biblabelprint{\csname\@citelabel{##1}\endcsname}%
         \@writeaux{\string\@citedef{##1}{\csname\@citelabel{##1}\endcsname}}%
         \ignorespaces
      }%
      %
      %
      \let\em = \bblem
      \let\newblock = \bblnewblock
      \let\sc = \bblsc
      \frenchspacing
      \clubpenalty = 4000 \widowpenalty = 4000
      \tolerance = 10000 \hfuzz = .5pt
      \everypar = {\hangindent = \biblabelwidth
                      \advance\hangindent by \biblabelextraspace}%
      \bblrm
      \parskip = 1.5ex plus .5ex minus .5ex
      \biblabelextraspace = .5em
      \bblhook
      \input \bblfilebasename.bbl
   \endgroup
}%
%
%
\@innernewdimen\biblabelwidth
\@innernewdimen\biblabelextraspace
%
%
%
\def\biblabelprint#1{%
   \noindent
   \hbox to \biblabelwidth{%
      \biblabelprecontents
      \biblabelcontents{#1}%
      \biblabelpostcontents
   }%
   \kern\biblabelextraspace
}%
%
%
%
\def\biblabelcontents#1{{\bblrm [#1]}}%
%
%
\def\bblrm{\rm}%
%
%
\def\bblem{\it}%
%
%
\def\bblsc{\ifx\@scfont\@undefined
              \font\@scfont = cmcsc10
           \fi
           \@scfont
}%
%
%
\def\bblnewblock{\hskip .11em plus .33em minus .07em }%
%
%
\let\bblhook = \empty
%
%
%
\def\printcitestart{[}
\def\printcitefinish{]}
\def\printbetweencitations{, }
\def\printcitenote#1{, #1}
%
%
%
\let\citation = \@gobble
%
%
%
\@innernewcount\@numparams
%
%
\def\newcommand#1{%
   \def\@commandname{#1}%
   \@getoptionalarg\@continuenewcommand
}%
%
%
\def\@continuenewcommand{%
   \@numparams = \ifx\@optionalarg\empty 0\else\@optionalarg \fi \relax
   \@newcommand
}%
%
%
\def\@newcommand#1{%
   \def\@startdef{\expandafter\edef\@commandname}%
   \ifnum\@numparams=0
      \let\@paramdef = \empty
   \else
      \ifnum\@numparams>9
         \errmessage{\the\@numparams\space is too many parameters}%
      \else
         \ifnum\@numparams<0
            \errmessage{\the\@numparams\space is too few parameters}%
         \else
            \edef\@paramdef{%
               \ifcase\@numparams
                  \empty  No arguments.
               \or ####1%
               \or ####1####2%
               \or ####1####2####3%
               \or ####1####2####3####4%
               \or ####1####2####3####4####5%
               \or ####1####2####3####4####5####6%
               \or ####1####2####3####4####5####6####7%
               \or ####1####2####3####4####5####6####7####8%
               \or ####1####2####3####4####5####6####7####8####9%
               \fi
            }%
         \fi
      \fi
   \fi
   \expandafter\@startdef\@paramdef{#1}%
}%
%
%
%
%
\def\@readauxfile{%
   \if@auxfiledone \else 
      \global\@auxfiledonetrue
      \@testfileexistence{aux}%
      \if@fileexists
         \begingroup
            \endlinechar = -1
            \catcode`@ = 11
            \input \jobname.aux
         \endgroup
      \else
         \message{\@undefinedmessage}%
         \global\@citewarningfalse
      \fi
      \immediate\openout\@auxfile = \jobname.aux
   \fi
}%
%
%
\newif\if@auxfiledone
\ifx\noauxfile\@undefined \else \@auxfiledonetrue\fi
%
%
%
%
\@innernewwrite\@auxfile
\def\@writeaux#1{\ifx\noauxfile\@undefined \write\@auxfile{#1}\fi}%
%
%
%
\ifx\@undefinedmessage\@undefined
   \def\@undefinedmessage{No .aux file; I won't give you warnings about
                          undefined citations.}%
\fi
%
%
\@innernewif\if@citewarning
\ifx\noauxfile\@undefined \@citewarningtrue\fi
%
%
%
\catcode`@ = \@oldatcatcode

\def\pfeilso{\leavevmode
            \vrule width 1pt height9pt depth 0pt\relax
           \vrule width 1pt height8.7pt depth 0pt\relax
           \vrule width 1pt height8.3pt depth 0pt\relax
           \vrule width 1pt height8.0pt depth 0pt\relax
           \vrule width 1pt height7.7pt depth 0pt\relax
            \vrule width 1pt height7.3pt depth 0pt\relax
            \vrule width 1pt height7.0pt depth 0pt\relax
            \vrule width 1pt height6.7pt depth 0pt\relax
            \vrule width 1pt height6.3pt depth 0pt\relax
            \vrule width 1pt height6.0pt depth 0pt\relax
            \vrule width 1pt height5.7pt depth 0pt\relax
            \vrule width 1pt height5.3pt depth 0pt\relax
            \vrule width 1pt height5.0pt depth 0pt\relax
            \vrule width 1pt height4.7pt depth 0pt\relax
            \vrule width 1pt height4.3pt depth 0pt\relax
            \vrule width 1pt height4.0pt depth 0pt\relax
            \vrule width 1pt height3.7pt depth 0pt\relax
            \vrule width 1pt height3.3pt depth 0pt\relax
            \vrule width 1pt height3.0pt depth 0pt\relax
            \vrule width 1pt height2.7pt depth 0pt\relax
            \vrule width 1pt height2.3pt depth 0pt\relax
            \vrule width 1pt height2.0pt depth 0pt\relax
            \vrule width 1pt height1.7pt depth 0pt\relax
            \vrule width 1pt height1.3pt depth 0pt\relax
            \vrule width 1pt height1.0pt depth 0pt\relax
            \vrule width 1pt height0.7pt depth 0pt\relax
            \vrule width 1pt height0.3pt depth 0pt\relax}

\def\pfeilsw{ \leavevmode 
            \vrule width 1pt height0.3pt depth 0pt\relax
            \vrule width 1pt height0.7pt depth 0pt\relax
            \vrule width 1pt height1.0pt depth 0pt\relax
            \vrule width 1pt height1.3pt depth 0pt\relax
            \vrule width 1pt height1.7pt depth 0pt\relax
            \vrule width 1pt height2.0pt depth 0pt\relax
            \vrule width 1pt height2.3pt depth 0pt\relax
            \vrule width 1pt height2.7pt depth 0pt\relax
            \vrule width 1pt height3.0pt depth 0pt\relax
            \vrule width 1pt height3.3pt depth 0pt\relax
            \vrule width 1pt height3.7pt depth 0pt\relax
            \vrule width 1pt height4.0pt depth 0pt\relax
            \vrule width 1pt height4.3pt depth 0pt\relax
            \vrule width 1pt height4.7pt depth 0pt\relax
            \vrule width 1pt height5.0pt depth 0pt\relax
            \vrule width 1pt height5.3pt depth 0pt\relax
            \vrule width 1pt height5.7pt depth 0pt\relax
            \vrule width 1pt height6.0pt depth 0pt\relax
            \vrule width 1pt height6.3pt depth 0pt\relax
            \vrule width 1pt height6.7pt depth 0pt\relax
            \vrule width 1pt height7.0pt depth 0pt\relax
            \vrule width 1pt height7.3pt depth 0pt\relax
            \vrule width 1pt height7.7pt depth 0pt\relax
            \vrule width 1pt height8.0pt depth 0pt\relax
            \vrule width 1pt height8.3pt depth 0pt\relax
            \vrule width 1pt height8.7pt depth 0pt\relax
            \vrule width 1pt height9pt depth 0pt\relax
      }


\def\widestnumber#1#2{}

\def\citewarning#1{\ifx\shlhetal\relax 
    \else
    \par{#1}\par
    \fi
}

\def\rm{\fam0 \tenrm}

\def\fakesubhead#1\endsubhead{\bigskip\noindent{\bf#1}\par}



%
%
%

%

\font\textrsfs=rsfs10
\font\scriptrsfs=rsfs7
\font\scriptscriptrsfs=rsfs5

\newfam\rsfsfam
\textfont\rsfsfam=\textrsfs
\scriptfont\rsfsfam=\scriptrsfs
\scriptscriptfont\rsfsfam=\scriptscriptrsfs

\edef\oldcatcodeofat{\the\catcode`\@}
\catcode`\@11

\def\Cal@@#1{\noaccents@ \fam \rsfsfam #1}

\catcode`\@\oldcatcodeofat


\expandafter\ifx \csname margininit\endcsname \relax\else\margininit\fi

\long\def\red#1\endred{}
\long\def\green#1\endgreen{}
\long\def\blue#1\endblue{}
\long\def\private#1\endprivate{}

\def\endred{ \unmatched endred! }
\def\endgreen{ \unmatched endgreen! }
\def\endblue{ \unmatched endblue! }
\def\endprivate{ \unmatched endprivate! }

\ifx\latexcolors\undefinedcs\def\latexcolors{}\fi

\def\emptycs{}
\def\evaluatelatexcolors{%
        \ifx\latexcolors\emptycs\else
        \expandafter\xxevaluate\latexcolors\xxfertig\evaluatelatexcolors\fi}
\def\xxevaluate#1,#2\xxfertig{\setupthiscolor{#1}%
        \def\latexcolors{#2}}


\font\smallfont=cmsl7
\def\rutgerscolor{\ifmmode\else\endgraf\fi\smallfont
\advance\leftskip0.5cm\relax}
\def\setupthiscolor#1{\edef\tmptmpcs{\noexpand\bgroup\noexpand\rutgerscolor
\noexpand\def\noexpand\currentcolor{#1}%
\noexpand}%
\expandafter\let\csname#1\endcsname\tmptmpcs
\def\tmptmpcs{\checkColorUnmatched{#1}\popthecolor}
\expandafter\let\csname end#1\endcsname\tmptmpcs}

\def\checkColorUnmatched#1{\def\expectcolor{#1}%
    \ifx\expectcolor\currentcolor   
    \else \edef\failhere{\noexpand\tryingToClose '\currentcolor' with end\expectcolor}\failhere\fi}

\def\currentcolor{???}

\def\popthecolor{\ifmmode\else\endgraf\fi\egroup}

\expandafter\def\csname#1\endcsname{}

\evaluatelatexcolors

 \let\outerhead\head
 \def\head{\innerhead}
 \let\innerhead\outerhead

 \let\outersubhead\subhead
 \def\subhead{\innersubhead}
 \let\innersubhead\outersubhead

 \let\outersubsubhead\subsubhead
 \def\subsubhead{\innersubsubhead}
 \let\innersubsubhead\outersubsubhead

 \let\outerproclaim\proclaim
 \def\proclaim{\innerproclaim}
 \let\innerproclaim\outerproclaim

 %
 %
 %
 %

\def\demo#1{\medskip\noindent{\it #1.\/}}
\def\enddemo{\smallskip}

\def\remark#1{\medskip\noindent{\it #1.\/}}
\def\endremark{\smallskip}

\pageheight{8.5truein}
\topmatter
\title{Model theoretic stability and categoricity for complete metric spaces}
 \endtitle
\rightheadtext{Model theoretic stability and categoricity}
\author {Saharon Shelah and Alexander Usvyatsov
\thanks {\null\newline We would like to thank 
 Alice Leonhardt for the beautiful typing. \null\newline
 The first author would like to thank the Israel Science Foundation for
 partial support of this research (Grant No. 242/03). 
 Publication 837. } \endthanks} \endauthor 

\affil{The Hebrew University of Jerusalem \\
Einstein Institute of Mathematics \\
Edmond J. Safra Campus, Givat Ram \\
Jerusalem 91904, Israel
 \medskip
 Department of Mathematics \\
 Hill Center-Busch Campus \\
  Rutgers, The State University of New Jersey \\
 110 Frelinghuysen Road \\
 Piscataway, NJ 08854-8019 USA
 \medskip
UCLA \\
 Department of Mathematics \\
P.O. Box 951555 \\
Los Angeles, CA 90095-1555 USA} \endaffil
\medskip

\abstract
We deal with the systematic development of stability for the context
of approximate elementary submodels of a monster metric space,
which is not far, but still very distinct from the first order
case.  In particular we 
prove the analogue of Morley's theorem for the classes of complete
metric spaces.  \endabstract
\endtopmatter
\document

\newpage

\head {\S1 Introduction and Preliminaries} \endhead  \resetall \sectno=1
 \spuriousreset
\bigskip

We work in the context of a compact homogeneous model ${\frak C}$ which is
also a complete metric space with a definable metric $\bold d(x,y)$
and all the predicates and function symbols respect the metric.  Such
a monster model will be called a ``monster metric space" (a momspace),
Definition \scite{MOM}.
We investigate the class $K$ of ``almost elementary" 
complete submodels of ${\frak C}$, Definition \scite{MOM.2.16}.  

This paper is devoted to categoricty of such classes
$K$ in uncountable cardinalities 
(generalizing Morley's theorem to this context).  
As we believe that isometry is too strong as a notion of
isomorphism for classes of metric structures, we try to weaken the
assumptions and work with $\varepsilon$-embeddings instead
(Definitions \scite{ISO.1}, \scite{ISO.3}). 

Several suggestions for a framework suitable for model theoretic
treatment of classes arising in functional analysis and dynamics have
been made in the last 40 years by Chang, Keisler, Stern, Henson, and
more recently by Ben-Yaacov.  All these attempts were concerned with
allowing a certain amount of compactness (e.g. ``capturing"
ultra-products of Banach spaces introduced by Krivine), without having
to deal with non-standard elements.  In this paper the authors choose
to work in the most general context which still allows compactness,
therefore generalizing all the above frameworks.  The main tools and
techniques used here are borrowed from homogeneous model theory.

Homogeneous model theory was introduced and first studied by Keisler,
developed further by the first author, Grossberg, Leslsmann, Hyttinen
and others.  It investigates classes of elementary (sometimes somewhat
``saturated") submodels of a big homogeneous model (``monster"), see
precise definitions later.  

In \cite{Sh:54}, the first author classified such ``monsters" with
respect to the amount of compactness they admit.  Monsters which are
compact in a language closed under negation are called ``monsters of
kind II".  Later Hrushovski suggested the name ``Robinson
Theories" for such classes, see \cite{Hrxz}.  Monsters which are compact
in a language not necessarily closed under negations are called
``monsters of kind III".  Recently Ben-Yaacov has studied this
context in great detail.  He called such monsters ``compact abstract
theories", in short CATS, see \cite{BY03}.

A simple generalization (replacing equality by definable metric)
allows us to speak of a monster model of a class of metric spaces.  We
call such monsters ``monster metric spaces", see Definition
\scite{MOM} below.  Several results are proven in this general
context, but some require compactness, Definition \scite{HMC.4}(2).
Therefore, our main theorem (Theorem \scite{CATEG}) holds in the
metric analogue of ``monsters of kind III", ``metric cats" in
Ben-Yaacov's terminology.

Independently of our work (and simultaneously), 
Ben-Yaacov investigated categoricity for
metric cats under the additional assumption of the topology on the
space of types being Hausdorff, which is the metric analogue of
``monsters of kind II" (Robinson theories).  In this context one
can reconstruct most of classical stability theory (e.g., independence
based on non-dividing), see \cite{BY05}.   These methods fail in the
more general context we were working at.  Therefore, techniques
developed and used here are very different, and rely heavily on
non-splitting and Ehrenfeucht-Mostowski constructions.  These tools do
not make any significant use of compactness, and we believe that this
assumption can be eliminated by modifying our methods slightly.
At this point we decided not to make the effort, but we try to mention
where exactly compactness is used.

In a recent work \cite{BeUs0y}, Ben-Yaacov and the second author
introduce a framework of continuous first order logic, closely
related to \cite{ChKe66} and show that once modified slightly, most
model-theoretic approaches to classes of metric spaces (such as
Henson's logic, see \cite{HeIo02}, Hausdorff metric cats, see
\cite{BY05}) are equivalent to continuous logic.  Although if continuous
model theory had been discovered earlier, this paper might have looked
differently, we would still like to point out that equivalences shown
in \cite{BeUs0y} do not include monster metric spaces, not even compact
ones.  The assumption of the topology on the type space being
Hausdorff is absolutely crucial in \cite{BeUs0y} and \cite{BY05}; it
provides us with the ability to ``approximate" negations, which makes
continuous logic very similar to classical first order logic (of
course, this has many advantages).  Lacking Hausdorff assumption, one 
has to use more general methods in order to reobtain basic properties.
This is why ``monsters of kind III" (general cats) have been studied
and understood much less than first order on Robinson theories, even
in the discrete (non-continuous, classical first order) context.  

Some work has been done, though: the first author proved the analogue
of Morley's theorem for existentially closed models in \cite{Sh:54},
classes of existentially closed models were investigated further by Pillay in
\cite{Pi00}, general cats were studied by Ben-Yaacov in \cite{BY03},
\cite{BY03a} and other papers.  
Our work continues this investigation in the more
general metric context (classical model theory can be viewed as a
particular case with discrete metric).

Several words should be said also about the difference between the
discrete and the metric context.  For example, why could we have not
simply modified the proof in \cite{Sh:54} slightly and obtain Theorem
\scite{CATEG}?  The answer is that our categoricity assumption is
significantly weaker, as we assume categoricity only for complete
structures.  For instance, the class of infinite-dimensional Hilbert
spaces is categorical in all densities, but not so is the class of
inner-product sapces, not necessarily complete.  Starting from a
weaker assumption we aim for a weaker conclusion; but consequences of
our assumption are still significantly weaker than in \cite{Sh:54}
(e.g. $\aleph_0$-stability is lost, we only have a topological
version), which complicates life significantly.

This work was originally carried out as a Ph.D. thesis of the second
author under the supervision of the first one.  The paper is an
 expanded version of the thesis wich was written in Hebrew and
submitted to the Hebrew University.
\bn
\centerline{$* \qquad * \qquad *$}
\bn

The paper is organized as follows:

We introduce our context in \S2.  In particular, we define the class
of models which is investigated (``almost elementary" submodels of
${\frak C},M \prec_1 {\frak C}$).  This notion generalizes Henson's
approximate elementary submodels (see \cite{HeIo02}).  It has the
following important property: if $M \prec_1 {\frak C}$, then its
completion (metric closure) $N$ also satisfies $N \prec_1 {\frak C}$.
We will be interested mostly in complete (as metric spaces) ``almost
elementary" submodels of ${\frak C}$.

The next section, \S3, is devoted to different kinds of approximations to
formulae and types.  The importance of considering these
approximations, i.e. topological neighborhoods of partial types, becomes clear
later, when stability, isolation, and other central notions are
discussed.

In \S4 we generalize the notion of stability in a cardinal $\lambda$
to our topological context.  We say that ${\frak C}$ is
$0^+-\lambda$-stable if for no $\varepsilon > 0$ can we find an
$\varepsilon$-net of $\lambda^+$ types over a set of cardinality
$\lambda$, i.e. the space of types over a set of cardinality $\lambda$
has (in a sense) density $\lambda$.  This is a generalization of
$\lambda$-stability for Banach spaces studied by Iovino in \cite{Io99}.
It is equivalent to the definitions given by Ben-Yaacov for Hausdorff
cats in \cite{BY05}, and in the context of continuous theories it 
coincides with the definitions given in \cite{BeUs0y}.

We prove equivalence of several similar definitions for
$0^+-\lambda$-stability, connect $0^+-\aleph_0$-stability to
non-splitting of types, classical stability in homogeneous model
theory, existence of average types, saturation of a (closure of a)
union of $(D,\aleph_1)$-homogeneous models.  Notions of isolation are
developed and density of strictly isolated types is proved under the
assumption of $0^+-\aleph_0$-stability.

In \S5 we develop the theory of Ehrenfeucht-Mostowski models in our
context.  As we lack forking calculus, some basic facts
(e.g. existence of $(D,\aleph_1)$-homogeneous models in all
uncountable cardinalities) require a different approach, which is
provided by the EM-models techniques.
Also in the proof of the main theorem (\S8) we take
advantage of the representation of the homogeneous model as an EM-model
in order to find inside it  a converging sequence which is ``close" to
a subsequence of a given uncountable sequence.

\S6 is devoted to notions of $\varepsilon$-embeddings and
$\varepsilon$-isomorphisms.  We try to weaken our assumptions as much as
we can, and choose to work with the following notion of ``weakly
uncountably categorical" classes: for every $\varepsilon > 0$ there
exists $\lambda > \aleph_0$ in which every two models and
$\varepsilon$-isomorphic to each other.  This property seems a priori
weaker than uncountable categoricity in some $\lambda > \aleph_0$,
and even than ``there exists $\lambda > \aleph_0$ such that for every
$\varepsilon > 0$ any two models of density $\lambda$ are
$\varepsilon$-isomorphic".  But the main theorem (\scite{CATEG})
states the following: assume ${\frak C}$ is weakly uncountably
categorical.  Then every model of density $\lambda > \aleph_0$ is
$(D,\lambda)$-homogeneous (in particular, unique up to an isometry).
So all the above notions turn out to be equivalent.

In \S7 we prepare the ground for the proof of \scite{CATEG}, showing that
any wu-categorical (weakly uncountably categorical) 
momspace is uni-dimensional (in the sense of
\cite{Sh:3}), i.e., any $(D,\aleph_1)$-homogeneous model of density
character $\lambda$ is, in fact, $(D,\lambda)$-homogeneous.  The last
section, \S8, contains the proof of the main theorem, Theorem \scite{CATEG}.

\bn
\centerline {$* \qquad * \qquad *$}
\bn
We recall now basic definitions considering homogeneous model theory
(see \cite{Sh:3},\cite{Sh:54} and \cite{GrLe02}):
\definition{\stag{0.1} Definition}  1) Let $\tau$ be a vocabulary,
$D$ a set of complete $\tau$-types over $\emptyset$ in finitely many
variables.  A $\tau$-model $M$ is called a \ub{$D$-model} if
$D(M) \subseteq D$ (where $D(M) =$ the set of complete
finite $\tau$-types over $\emptyset$ realized in $M$).  $M$ is called
\ub{$(D,\lambda)$-homogeneous} if $D(M) = D$ and
$M$ is $\lambda$-homogeneous.
\nl
2) We call $D$ as in (1) a \ub{finite diagram} if for some model
$M,D = D(M)$.
\nl
3) $A$ is a $D$-set in $M$ if $A \subseteq M$ and $\bar a \in
{}^{\omega >}A \Rightarrow \text{\rm tp}(\bar a,\emptyset,M) \in D$.
For $A$ a $D$-set, $p \in \bold S^m_D(A,M)$ if there are $N,\bar a$ such
that $M \prec N,\bar a \in N$ realizes $p$ and $A \cup \bar a$ is a
$D$-set in $N$.
\nl
\enddefinition
\bigskip

\remark{\stag{0.2} Remark}  Let $D$ be a finite diagram.  $M$ is
$(D,\lambda)$-homogeneous \ub{iff} $M$ is universal for 
$D$-models of cardinality $\le \lambda$ (or just $\lambda$-sets when
$\lambda < |\tau| + \aleph_0$)  and $\lambda$-homogeneous
\ub{iff} $D(M)=D$ and $M$ is 
$(D,\lambda)$-saturated (i.e. every $D$-type over a 
subset of cardinality $< \lambda$ is realized in $M$).
\endremark
\bigskip

\demo{Proof}  Easy (or see \cite{GrLe02}(2.3),(2.4)).
\enddemo
\bn
This motivates the following
\definition{\stag{0.3} Definition}  Let $\lambda^*$ be big enough.  A
$(D,\lambda^*)$-homogeneous model ${\frak C}$ will be called a
\ub{$D$-monster} (or a homomogeneous monster for $D$) .  
We usually assume $\|{\frak C}\| = \lambda^*$.
\enddefinition
\bn
Recall:
\definition{\stag{0.4} Definition}  1) A finite diagram $D$ is
called $\lambda$-\ub{good} if there is a
$(D,\lambda)$-homogeneous model $M$ of cardinality $\ge \lambda$.
\nl
2) $D$ is good if it is $\lambda$-good for every $\lambda$.
\enddefinition
\bigskip

\remark{\stag{0.5} Remark}  So $D$ is good iff it has a monster.
\endremark
\bigskip

\demo{\stag{0.6} Convention}  In this paper we will fix a 
good finite diagram $D$ and a $D$-monster model ${\frak C}$.
\enddemo
\bn
\bigskip

\demo{\stag{0.7} Observation}  1) If ($D$ is good, ${\frak C}$ a
$D$-monster model) $p \in \bold S_D(A),A \subseteq B$ \ub{then}
there is $q,p \subseteq q \in \bold S_D(B)$.
\enddemo
\bn
\ub{Question}:  Why can we use ``good $D$"?  There are several
answers.

Basically, Claim \scite{0.7S} says that every $D$ which is the
finite diagram of a \ub{compact} momspace (see Definitions
\scite{MOM}, \scite{MOM.1}) is good.  Claim \scite{0.8S} says that
even without compactness, categoricity implies stability, which
implies $D$ is good by \cite{Sh:3}.  The reader can easily omit these
claims in the first reading.
\bigskip

\proclaim{\stag{0.7S} Claim}  Assume
\mr
\item "{$(a)$}"  $D$ is a finite diagram, i.e., a set of complete
$n$-types in $\Bbb L(\tau_{\frak C})$
\sn
\item "{$(b)$}"  ${\frak C}$ is $(D,\aleph_0)$-homogeneous
\sn
\item "{$(c)$}"  $\Delta$ is full for ${\frak C}$ (see Definition
\scite{HMC.4}) 
\sn
\item "{$(d)$}"  ${\frak C}$ is $\Delta$-compact.
\ermn
\ub{Then}
\mr
\item "{$(*)$}"  $(\alpha) \quad D$ is good
\sn
\item "{${{}}$}"  $(\beta) \quad$ if $\kappa(D) < \infty$ then $\kappa(D) \le
(|\tau_{\frak C}| + \aleph_0)^+$ ($\kappa(D)$ as in \cite{Sh:54}).
\endroster
\endproclaim
\bigskip

\demo{Proof}  See \cite{Sh:54}, the discussion of monsters of kind III or
\cite{BY03}, existence of a universal domain.
\enddemo
\bigskip

\proclaim{\stag{0.8S} Claim}  1) Assume that
\mr
\item "{$(a)$}"  $\tau$ a countable metric vocabulary 
and let $\delta(*) = (2^{\aleph_0})^+$
\sn
\item "{$(b)$}"  $\tau = \tau_D$
\sn
\item "{$(c)$}"  $D$ is a finite diagram
\sn
\item "{$(d)$}"  there is a $(D,\aleph_1)$-homogeneous model $M$
\sn
\item "{$(e)$}"  there is a $D$-model of cardinality $\ge 
\beth_{\delta(*)}$ (probably less is enough) and
\sn
\item "{$(f)$}"   $K^1_D = \{M:M \in K^c_1\}$ (see Definition \scite{MOM.2.16})
 is categorical in $\lambda > \aleph_0$ 
or is wu-categorical (see Definition \scite{ISO.3}).
\ermn
\ub{Then} $D$ is stable, hence is good. 
\endproclaim
\bigskip

\demo{Proof}  By \scite{OS}, \scite{STAB.2s.0} hence we get stability,
$D$ is good follows by \cite{Sh:3}.
\enddemo
\bigskip

\demo{\stag{0.9} Notations}  

$\lambda,\mu,\chi \qquad \qquad$ infinite cardinals

$\alpha,\beta,\gamma \qquad \qquad$  infinite ordinals

$\delta \qquad \qquad \qquad$ limit ordinals

$\nu,\eta \qquad \qquad \qquad$ sequences of ordinals

$\varphi,\psi,\vartheta \qquad \qquad$ formulae

${\frak C} \qquad \qquad \quad$ the monster model

$M,N \qquad \qquad$ models (in the monster)

$A,B,C \qquad \qquad$ sets (in the monster)

$\varepsilon,\zeta,\xi \qquad \qquad$ positive reals

$I,J \qquad \qquad$ order types

$\bold I,\bold J \qquad \qquad$ indiscernible sequences

$\bold d \qquad \qquad \qquad$ a metric
\enddemo
\bigskip

\head {\S2 Main context - monster metric spaces} \endhead  \resetall \sectno=2
 \spuriousreset
\bigskip

In this section we discuss our main context.  We start with some
notations.
\definition{\stag{1n.1} Definition}  Let $(X,\bold d)$ be a 
metric space.
\nl
1) We extend the metric to
$n$-tuples: for $\langle a_i:i < n \rangle,\langle b_i:i < n\rangle$
we define $\bold d(\bar a,\bar b) = \text{\rm max}\{\bold
d(a_i,b_i):i<n\}$.
\nl
2) For $\bar a$ an $n$-tuple, $A$ a set of $n$-tuples, 
$\bold d(\bar a,A) = \bold d(A,\bar a) = \text{\rm inf}\{\bold d(\bar
a,\bar b):\bar b \in A\}$.
\nl
3) For sets $A,B \subseteq X^n$, we define two versions
of distances:

$$
\bold d_1(A,B) = \text{\rm inf}\{\bold d(\bar a,\bar b):\bar a \in
A,\bar b \in B\}
$$

$$
\bold d_2(A,B) = \text{\rm sup}\{\{\bold d_1(\bar a,B):\bar a \in A\}
\cup \{\bold d_1(A,\bar b):\bar b \in B\}\}.
$$
\mn
4) We denote the density (the density character) of $X$ by Ch$(X)$.
So Ch$(X)$ is the minimal cardinality of a dense subset of $X$.
\nl
5) We denote the topological (metric) closure of a set $A$ by $\bar A$
or mcl$(A)$.
\enddefinition
\bigskip

\definition{\stag{1.1} Definition}  1) We call a vocabulary $\tau$
\ub{metric} if it contains predicates $P_{q_1,q_2}(x,y)$ for all $0
\le q_1 \le q_2$ rationals.  We call the collection of these
predicates a \ub{metric scheme}.
\nl
2) Given a metric vocabulary $\tau$, we call a $\tau$-structure 
$M$ \ub{semi-metric} if for some (unique) $\bold d$:
\mr
\item "{$(i)$}"    $M$ is a metric space with the metric $\bold d$
\sn
\item "{$(ii)$}"  $\bold d$ is definable in $M$ by the $\tau$-metric
scheme, i.e., $\bold d(a,b) \in [q_1,q_2]$ iff $M \models
P_{q_1,q_2}(a,b)$ for all rationals $q_1,q_2$ and $a,b \in M$.
\ermn
3) We call a semi-metric structure \ub{complete} if $(M,\bold d)$ is
complete as a metric space.
\enddefinition
\bigskip

\remark{\stag{1.2} Remark}  Given a semi-metric $\tau$-structure $M$, we
will usually write $M \models ``\bold d(a,b) \le q"$, etc., forgetting
the $\tau$-metric scheme.

Note that in a pseudo-metric structure for each $r_1,r_2 \in \Bbb R$, the
property ``$r_1 \le \bold d(x,z) \le r_2$" is $0$-type-definable by a
set of atomic formulas.
\endremark
\bigskip

\definition{\stag{HMC.sO} Definition}  1) $(M,\bold d)$ is a 
\ub{metric} structure (or model) \ub{when}:
\mr
\item "{$(a)$}"  $\tau({\frak C})$ is a metric vocabulary and
$(M,\bold d)$ is a semi-metric structure
\sn
\item "{$(b)$}"  $P^M \subseteq {}^{\text{arity}(P)}M$ is closed
(with respect to $\bold d$) for every predicate $P \in \tau_M$
\sn
\item "{$(c)$}"  $F^M:{}^{\text{arity}(F)}M \rightarrow N$ is a
continuous function for every function symbol $F \in \tau_M$.
\ermn
2) $({\frak C},\bold d)$ is a (homogeneous) metric monster 
if: $({\frak C},\bold d)$
is a metric model, ${\frak C}$ is a $D$-monster for some finite
diagram $D$,  constant here.
\enddefinition
\bigskip

\remark{\stag{1.3} Remark}  Each homogeneous monster admits the
discrete definable metric, i.e., $\bold d(a,b)=1$ for all $a \ne b \in
{\frak C}$, and it is definable by the equality and inequality.  So
each homogeneous monster is a metric monster with the discrete metric.
\endremark
\bigskip

\demo{\stag{1.4} Notations}  We will often identify $\varphi = 
\varphi(\bar x,\bar a)$ (for $\varphi$ a formula, 
$\bar a \in {\frak C}$) with the set
of the realizations in ${\frak C}$ of $\varphi(\bar x,\bar a)$, i.e.,
$\varphi({\frak C},\bar a) = \varphi^{\frak C} =: \{\bar b \in {}^{\ell g(\bar
x)}{\frak C}:{\frak C} \models \varphi[\bar b,\bar a]\}$.
So $\bold d_1(\varphi,\psi)$ in
fact means $\bold d_1(\varphi^{\frak C},\psi^{\frak C})$, etc.
\enddemo
\bn
\centerline {$* \qquad * \qquad *$}
\bn
The following definition is the analogue of abstract elementary classes
(see \cite{Sh:88r}) in our context.
\definition{\stag{AMC} Definition}  Let $({\frak K},\le_{\frak K})$ be
an ordered class of $\tau = \tau({\frak K})$ complete metric
structures ($\tau$ is a metric vocabulary), ${\frak K}$ closed under
$\tau$-isomorhisms.  We call \nl
$({\frak K},\le_{\frak K})$ an
\ub{abstract metric class} (a.m.c.) if
\mr
\item   $\langle M_i:i < \delta \rangle$ is a 
$\le_{\frak K}$-increasing sequence (from ${\frak K}$, of
course), then $M = \text{\rm mcl}(\dbcu_{i<\delta} M_i) 
\in {\frak K}$.  Moreover, $[i < \delta \Rightarrow M_i \le_{\frak K}
M]$ and $[M_i \le_{\frak K} N$ for all $i<\delta 
\Rightarrow M \le_{\frak K} N]$
\sn
\item  for $M_1 \subseteq M_2 \subseteq M_3$ from ${\frak K},[M_1,M_2
\le_{\frak K} M_3] \Rightarrow M_1 \le_{\frak K} M_2$
\sn
\item  for some cardinal LS$({\frak K})$ we have the ``downward
L\"owenheim-Skolem theorem", i.e., for each $M \in {\frak K},A
\subseteq M$, there exists $N \in {\frak K},N \le_{\frak K} M,A
\subseteq N$, Ch$(N) \le |A| + \text{\rm LS}({\frak K})$.
\endroster
\enddefinition
\bigskip

\definition{\stag{HMC.2} Definition}  We say that a set of formula
$\Delta$ for a homogeneous metric monster ${\frak C},\bold d$ is
\ub{admissible} if for each $\varphi(\bar x) \in \Delta$, the set
$\varphi^{\frak C} = \{\bar a \in {\frak C},{\frak C} \models
\varphi(\bar a)\}$ is closed with respect to the metric $\bold d$
(topology induced by it). 
\enddefinition
\bigskip

\definition{\stag{HMC.2S} Definition}  In ${\frak C},\varphi,\psi$ are
contradictory if $\bold d_1(\varphi^{\frak C},\psi^{\frak C}) > 0$.
\enddefinition
\bigskip

\remark{\stag{HMC.2T} Remark}  Later (see \scite{MOM.3A}) we show that
in our context this is equivalent to $\varphi^{\frak C} \cap
\psi^{\frak C} = \emptyset$.
\endremark
\bigskip

\demo{\stag{HMC.3} Example}  If $\bold d$ is discrete, then each subset
of ${\frak C}$ is closed, so the set of all formulas $\Delta = \Bbb
L(\tau_{\frak C})$ is admissible.
\enddemo
\bn
In order to give a non-trivial example, we define
\definition{\stag{POS.1} Definition}  For a metric model $(M,\bold d)$
and a formula $\varphi(x,\bar a)$ with parameters $\bar a \in M$, we
say that $M \models \exists^* x \varphi(x,\bar a)$ (there almost
exists $x$ such that $\varphi(x,\bar a)$) if for every $\varepsilon > 0$ there
exist $b,\bar a'$ such that $M \models \varphi(b,\bar a')$ and $\bold
d(\bar a,\bar a') \le \varepsilon$.
\enddefinition
\bigskip

\definition{\stag{POS} Definition}  1) We define \ub{positive} formulae
by induction: each atomic formula is positive, for $\varphi,\psi$
positive, $\varphi \wedge \psi,\varphi \vee \psi,\exists^* x
\varphi,\forall x \varphi$ are positive.  So negation and implication
are not allowed.
\nl
2) \ub{Positive existential} formulae are defined  similarly without
allowing $\forall x \varphi$.
\enddefinition
\bigskip

\demo{\stag{HMC.so.2} Observation}  For $(M,\bold d)$ a metric model
and positive $\varphi(\bar x) \in \Bbb L(\tau_M),\varphi^M = \{\bar a
\in {}^{\ell g(\bar x)}M:M \models \varphi[\bar a]\}$ is a closed
subset of ${}^{\ell g(\bar x)}M$ (under $\bold d$), i.e., $\Delta =$
``the positive formulae" is admissible.
\enddemo
\bigskip

\definition{\stag{HMC.4} Definition}  1) For a homogeneous monster
${\frak C}$, a set of formulae $\Delta$ is called \ub{full} if
\mr
\item "{$(i)$}"  for all $\bar a \in {\frak C},A \subseteq {\frak C}$,
tp$(\bar a,A,{\frak C})$ is determined by the $\Delta$-type
tp$_\Delta(a,A,{\frak C})$
\sn
\item "{$(ii)$}"  if $\varphi \in \Delta,\bar a \in {\frak C}$ and
${\frak C} \models \neg \varphi(\bar a)$, then there exists $\psi \in
\Delta$ such that $(\varphi,\psi)$ is contradictory (see \scite{HMC.2S})
and ${\frak C} \models \psi(\bar a)$.
\ermn
2) We call ${\frak C}$ \ub{$\Delta$-compact} (where $\Delta$ is a set
of formulae, i.e., $\subseteq \Bbb L(\tau_{\frak C})$) 
if each set of $\Delta$-formulae with parameters from
${\frak C}$ of cardinality $<|{\frak C}|$ which is finitely
satisfiable in ${\frak C}$, is realized in ${\frak C}$.  We omit
$\Delta$ if constant, and abusing notation write $\Delta = \Delta({\frak C})$.
\enddefinition
\bigskip

\definition{\stag{HMC.4s.0} Definition}  We say that $\Delta$ is
full$^+$ if: as in \scite{HMC.4} but in (ii) the quantifier depth of
$\psi$ is $\le$ the quantifier depth of $\varphi$.
\enddefinition
\bn
Now we make the main definition of this section, introducing the
context of this paper.
\definition{\stag{MOM} Definition}  1) A metric 
homogeneous class $({\frak K},
\le_{\frak K})$ (equivalently: its metric homogeneous monster
${\frak C}$) is called $\Delta$-\ub{momspace} (monster metric space)
if $\Delta = \Delta({\frak C})$ is a set of formulae 
containing the metric scheme, such that
\mr
\item "{$(a)$}"    $\Delta \subseteq \Bbb L(\tau_{\frak C})$ is closed
under $\wedge,\exists^*$ and subformulae 
\sn
\item "{$(b)$}"  $\Delta$ is admissible
\sn
\item "{$(c)$}"  $\Delta$ is full for ${\frak C}$.
\ermn
2) Let ``${\frak C}$ is momspace" mean ``for some $\Delta$" and choose
such $\Delta = \Delta({\frak C})$ (well, it is not necessarily unique
but we ignore this).
\nl
3) The metric $\bold d = \bold d_{\frak C}$ can be defined from
${\frak C}$ so we can ``forget" to mention $\bold d$, but still will
usually say $({\frak C},\bold d)$, e.g. to distinguish from 
${\frak C}$ when we use the
$(D,\lambda)$-homogeneous context.
\enddefinition
\bigskip

\demo{\stag{MOM.2.15} Convention}  $({\frak C},\bold d)$ is a fixed momspace.
\enddemo
\bn
The class of models we are interested in is defined below.  
\definition{\stag{MOM.2.16} Definition}   1) $K_1 = K^1_{\frak C}$ is
the class of $M$ such that:
\mr
\item "{$(a)$}"  $M \subseteq {\frak C}$ 
\sn
\item "{$(b)$}"  $M \prec^1_{\Delta} {\frak C}$ which means:
\nl
\ub{if} $\varepsilon > 0,{\frak C} \models \varphi[\bar b,\bar
a],\varphi(\bar y,\bar x) \in \Delta,\bar a \in {}^{\ell g(\bar x)}M$
\ub{then} for some $\bar b' \in {}^{\ell g(\bar y)} M$ we have 
${\frak C} \models (\exists \bar x,\bar y)(\varphi(\bar y,\bar x) \wedge
\bold d(\bar y \char 94 \bar x,\bar b' \char 94 \bar a)$.
\ermn
2) $K^c_1$ is the class of members of $K_1$ which are \ub{complete}.
\enddefinition
\bigskip

\proclaim{\stag{MOM.2.19} Claim}  1) Assume $M \subseteq {\frak C}$ and
$N = \text{\rm mcl}(M)$.  Then
\mr
\item "{$(a)$}"  $M \in K_1$ iff $N \in K_1$ iff $N \in K^c_1$.
\ermn
2) $(K_1,\subseteq)$ is an abstract metric class 
\endproclaim
\bigskip

\definition{\stag{MOM.1} Definition}  A momspace $({\frak C},\bold d)$
is called \ub{compact} if ${\frak C}$ is $\Delta({\frak C})$-compact.
\enddefinition
\bigskip

\demo{\stag{MON.1s.2} Convention}   We work in a compact
$\Delta$-momspace $({\frak C},\bold d),M,N$ denote submodels which
are from $K_1$; though really interested in the closed (complete) ones,
i.e. versions of categoricity are defined using complete models.
\ub{BUT} we try to mention when we use compactness.
\enddemo
\bigskip

\demo{\stag{MOM.1s.3} Observation}  [${\frak C}$ compact]  For
$\varphi(x,\bar y) \in \Delta,\bar a \in {\frak C}$, we have ${\frak
C} \models \exists^* x \varphi(x,\bar a) \Leftrightarrow {\frak C}
\models \exists x \varphi(x,\bar a)$ [so we can ``forget" about the
new existential quantifier and use the classical one].
\enddemo
\bigskip

\demo{Proof}  Assume ${\frak C} \models \exists^* x \varphi(x,\bar
a)$.  Then the set $\{\varphi(x,\bar y),\bold d(\bar y,\bar a) \le
\frac 1n:n < \omega\}$ is consistent in ${\frak C}$ (by compactness),
so ${\frak C} \models \exists x \varphi(x,\bar a)$. \hfill$\square$
\enddemo
\bigskip

\demo{\stag{MOM.1.2.19} Examples}  1) Let $T$ be first order, ${\frak C}$ a
big saturated model of $T$, then ${\frak C}$ is a compact momspace
(with discrete metric), $\Delta({\frak C}) = {\Cal L}$ (all formulae).
\nl
2) Let $T$ be a Robinson theory (see \cite{Hrxz}), ${\frak C}$ its
universal domain.  Then it is a compact momspace.
\nl
3) Consider the unit ball of a monster Banach space as in 
\cite{HeIo02}, Chapter 12.
Then $({\frak C},\|\cdot\|$) is a compact momspace (we write the norm
here instead of the metric) where $\Delta({\frak C}) =$ positive formulae (more
precisely positive bounded formulae, see \cite{HeIo02}, Chapter 5),
${\frak K} = K^c_1 = \{M:M \prec_A {\frak C}\}$, see
\cite{HeIo02}(6.2) for the definition of $\prec_A$.
Compactness is proven in Chapter 9 there.  One can easily show that
$\prec_A = \prec^1_\Delta$, see also \scite{MOM.suB}, \cite{HeIo02}(6.6).
\nl
4) Metric Hausdorff CATs, see \cite{BY03} for a definition of a Hausdorff
 CAT.  A CAT is metric if its monster is metric.  See \cite{BY05} in
Hausdorff metric cats.
\nl
5) The main example we have in mind: $({\frak C},\bold d)$ is a
compact momspace for $\Delta({\frak C}) =$ positive existential
formulae, ${\frak K} = K^c_1$.
So in particular it is a metric cat (not necessarily Hausdorff).
\nl
\enddemo
\bigskip

\head {\S3 Approximations to formulas and types} \endhead  \resetall \sectno=3
 \spuriousreset
\bigskip

The following notions of $\varepsilon$-approximations of formulae is of
major importance.  We give two different definitions and will use both
for different purposes.  Note that \scite{MOM.2} simply defines topological
neighborhoods, while in \scite{MOM.2s.13} by moving the parameters we 
allow the formula to change
a little, see also \scite{MOM.suB}.  Let $\varepsilon$ denote a
non-negative rational number.
\bigskip

\definition{\stag{MOM.2} Definition}  1) For a formula $\varphi(\bar
x)$ possibly with parameters, we define 
$\varphi^{[\varepsilon]}(\bar x) = \exists \bar x'(\varphi(\bar x') \wedge
\bold d(\bar x,\bar x') \le \varepsilon)$.  So
$\varphi^{[\varepsilon]}(\bar x,\bar a) = (\exists \bar x')[\bold
d(\bar x,\bar x') \le \varepsilon \wedge \varphi(\bar x',\bar a)]$. 
\nl
2) For a (partial) type $p$, define $p^{[\varepsilon]} = 
\{(\dsize \bigwedge_{\ell < n} \varphi_\ell)^{[\varepsilon]}:n < \omega$ and
$\varphi_0,\dotsc,\varphi_{n-1} \in p\}$.
\nl
\enddefinition
\bigskip

\definition{\stag{MOM.2s.13} Definition}  1) For a formula 
$\varphi(\bar x,\bar a)$ let $\varphi^{<\varepsilon>}(\bar
x,\bar a) = (\exists \bar x',\bar y')(\varphi(\bar x',\bar y')
\wedge \bold d(\bar x',\bar x) \le \varepsilon \wedge 
\bold d(\bar y',\bar a)  \le \varepsilon)$.
\nl
2) For a (partial) type $p,p^{<\varepsilon>}(\bar x) = \{(\dsize \bigwedge_\ell
\varphi_\ell(\bar x,\bar a_\ell))^{<\varepsilon>}:n < \omega$ and
$\varphi_\ell(\bar x,\bar a_\ell) \in p$ for $\ell <n\}$. 
\enddefinition
\bigskip

\definition{\stag{MOM.2T} Definition}  1) For a (partial) type (maybe
with parameters) $p$, or set $B$ and $\varepsilon > 0$, we say that
$\bar c \in B$ realizes $p^{[\varepsilon]}$ if ${\frak C} \models
p^{[\varepsilon]}(\bar c)$.  We say that $p^{[\varepsilon]}$ is
realized in $B$ if some $\bar c \in B$ realizes it.
\nl
2) Similar for $p^{<\varepsilon>}$. 
\enddefinition
\bigskip

\remark{\stag{MOM.2s.3} Remark}  1) Note that
$\varphi^{[\varepsilon]}(\bar x,\bar a)$ is equivalent to
$\varphi^{[\varepsilon]}_1(\bar x)$, when we expand ${\frak C}$ by
individual constants for $\bar a$ and let $\varphi_1(\bar x) 
= \varphi(\bar x,\bar a)$.
\nl
2) $\varphi^{[\varepsilon]}\models \varphi^{<\varepsilon>}$ for all 
$\varphi$.
\nl
3) For a formula without parameters $\varphi,\varphi^{[\varepsilon]} =
\varphi^{<\varepsilon>}$.
\nl
4) For $\zeta \ge \varepsilon \ge 0,\varphi^{[\varepsilon]} \models
\varphi^{[\zeta]},\varphi^{<\varepsilon>} \models
   \varphi^{<\zeta>}$.
\nl
5) For $\varepsilon=0,\varphi^{[\varepsilon]} \equiv \varphi^{<\varepsilon>}
\equiv \varphi$.
\endremark
\bigskip

\demo{\stag{MOM.3} Observation}  1) If $\varphi = \varphi(\bar x,\bar b)$ 
(as usual is admissible), then
$\varphi \equiv \dsize \bigwedge_{\varepsilon >0}
\varphi^{[\varepsilon]}$.
\nl
2) If $\varphi = \varphi(\bar x,\bar a)$ then $\varphi  
\equiv \dsize \bigwedge_{\varepsilon > 0} \,\,
\dsize \bigvee_{\zeta \in (0,\varepsilon)} \varphi^{[\zeta]}$.
\nl
3) Similar for $\varphi^{<\varepsilon>}$.
\enddemo
\bigskip

\demo{Proof}  1) $[a \models \varphi \Rightarrow a \models \bigwedge
\{\varphi^{[\varepsilon]}:\varepsilon > 0\}]$ is obvious.

Suppose $\bar a \models \varphi^{[\varepsilon]}$ for all $\varepsilon$, so
for each $n$ there exists $\bar a_n$ such that ${\frak C} 
\models \varphi[\bar a_n]$ and $\bold d(\bar a,\bar a_n) 
\le \frac 1n$.  So $\langle \bar a_n:n \in \omega 
\rangle$ converges to $\bar a$,
therefore $\bar a \models \varphi$, as $\varphi^{\frak C}$ is a closed
set.
\nl
2),3)  Similar.
\enddemo
\bigskip

\proclaim{\stag{4n.10} Claim}   1) $(\theta^{<\zeta_1>}(\bar
x))^{<\zeta_2>} \equiv \theta^{<\zeta_1+\zeta_2>}(\bar x)$.
\nl
2) $(p^{<\zeta_1>}(\bar x))^{<\zeta_2>} \equiv 
p^{<\zeta_1+\zeta_2>}$.
\nl
3) Assume $\bold d(\bar b,\bar b') \le \zeta_1 - \zeta_2$ then
   $\theta^{<\zeta_2>}(x,b') \models \theta^{<\zeta_1>}(x,b)$.
\endproclaim
\bigskip

\proclaim{\stag{MOM.3A} Claim}  1) If ${\frak C}$ is $\Delta$-compact
\ub{then} a pair (of admissible $\Delta$-formulas) $\varphi 
= \varphi(\bar x,\bar a),\psi 
= \psi(\bar x,\bar b)$ is contradictory (see definition \scite{HMC.2S})
iff $\varphi^{\frak C} \cap \psi^{\frak C} = \emptyset$. 
\nl
2) [${\frak C}$ not necessarily compact]  A pair as in (1)
$\varphi,\psi$ is contradictory iff for some $\varepsilon >0,
\varphi^{[\varepsilon]} \cap \psi^{[\varepsilon]} = \emptyset$.
\nl
3) If $\varphi,\psi$ in (2) are without parameters, we can add ``iff
   $\varphi^{<\varepsilon>} \cap \psi^{<\varepsilon>} = \emptyset$ for
   some $\varepsilon > 0$".
\nl
4) If ${\frak C}$ is compact, (3) holds for formulae with parameters.
\endproclaim
\bigskip

\demo{Proof}  1) Obviously if $(\varphi,\psi)$ is a contradictory pair
then $\varphi^{\frak C} \cap \psi^{\frak C} = \emptyset$.  Now assume
$\bold d_1(\varphi,\psi) = 0$.  So for each $n > 0$ there exists
$\bar a_n,\bar b_n \in {\frak C}$ such that $\bold d(\bar a_n,\bar
b_n) < \frac 1n$ and
$\varphi(\bar a_n,\bar a),\psi(\bar b_n,\bar b)$ hold.  Therefore the set
$\{\varphi(\bar x,\bar a),\psi(\bar y,\bar b)\} \cup \{\bold d(\bar
x,\bar y) \le \frac 1n:n=1,2,\ldots\}$ is finitely consistent.  By
compactness we obtain $\bar a^*,\bar b^*$ realizing it, but
necessarily $\bold d(\bar a^*,\bar b^*) = 0$, so $\bar a^* = \bar b^*$
realizes both $\varphi$ and $\psi$ and we are done.
\nl
2) Trivial.
\nl
3) Follows from (2) and \scite{MOM.2s.3}(3).
\nl 
4) Assume $\varphi^{<\varepsilon>} \cap \psi^{<\varepsilon>} \ne
\emptyset$ for all $\varepsilon > 0$, so the set
$\{\varphi^{<\varepsilon>}(\bar x,\bar a) \wedge
\psi^{<\varepsilon>}(\bar x,\bar b):\varepsilon >0\}$ is finitely
satisfiable (by \scite{MOM.2s.3}(4)) and therefore realized in ${\frak
C}$, now by \scite{MOM.3}(3) $\varphi \cap \psi \ne \emptyset$,
therefore $\varphi$ and $\psi$ are not contradictory.  The other
direction is trivial.  \hfill$\square_{\scite{MOM.3A}}$
\enddemo
\bigskip

\proclaim{\stag{MOM.4} Claim}  [Let $({\frak C},\bold d)$ be a compact
momspace], $p$ a type over a set $A$.
\nl
Then $\bar a \in {\frak C}$ realizes $p^{[\varepsilon]}$ iff there
exists $\bar a' \in {\frak C}$ realizing $p,\bold d(\bar a,\bar a')
\le \varepsilon$ (so $p^{[\varepsilon]}$ is just the
``$\varepsilon$-neighborhood" of $p$).
\endproclaim
\bigskip

\demo{Proof}  The if direction is obvious.  Now suppose $\bar a
\models p^{[\varepsilon]}$, so $[\varphi(\bar a)]^{[\varepsilon]}$
holds for all $\varphi \in p$, but there is no $\bar a' \models p$
such that $\bold d(\bar a,\bar a') \le \varepsilon$.  Then the set

$$
\{\bold d(\bar x',a') \le \varepsilon\} \cup \{\varphi(\bar x'):
\varphi \in p,\varphi \text{ is a } \Delta({\frak C})\text{-formula}\}
$$
\mn
is inconsistent (as $\Delta$ is full for ${\frak C}$), and by
compactness we get a contradiction (because in the definition of
$p^{[\varepsilon]}$ we essentially close $p$ under conjunctions).
\enddemo
\bn
We connect the relevant notion of submodel defined in \scite{MOM.2.16}
with the ``logical" notion of approximation:
\demo{\stag{MOM.suB} Observation}  $M \prec^1_\Delta N$ ($M,N$
structures, see \scite{MOM.2.16} for the definition of
$\prec^1_\Delta$) iff for every $\varphi(\bar y,\bar x) \in
\Delta,\varepsilon > 0$, if $N \models \exists \bar y \varphi(\bar
y,\bar a),\bar a \in M$, then $N \models \varphi^{<\varepsilon>}(\bar
b,\bar a)$ for some $\bar b \in M$.
\enddemo
\bigskip

Note that the assumption above of $M$ being complete in fact 
follows from $(D,\aleph_1)$-homogeneous:
\bigskip

\demo{\stag{MOM.sub.2} Observation}  Let $M \in K_1$ be
$(D,\aleph_1)$-homogeneous, then $M \in K^c_1$, i.e., is complete.
\enddemo
\bigskip

\demo{Proof}  Let $\langle a_n:n < \omega\rangle$ 
converge to $a \in {\frak C},a_n \in M$.  So the set 
$\{\bold d(x,a_n) \le {\bold d} (a,a_n):n \in \Bbb
N\}$ is realized in $M$, so $\langle a_n:n <\omega\rangle$ converges
in $M$ (to the same limit, of course).  \hfill$\square_{\scite{MOM.sub.2}}$
\enddemo
\bigskip

\definition{\stag{MOM.4A} Definition}  We call $M \in K^c_1$ pseudo
$(D,\lambda)$-homogeneous if for every $A \subseteq M,|A| < \lambda$,
for every $p \in S^m_D(A)$, for every $\varepsilon  >
0,p^{<\varepsilon >}(\lambda)$ is realized in $M$ (see Definition
\scite{MOM.2T}).
\nl
Note that replacing $p^{<\varepsilon>}$ above by $p^{[\varepsilon]}$,
we get $(D,\lambda)$-homogeneous, see \scite{MOM.4B} below. 
\enddefinition
\bigskip

\definition{\stag{MOM.4A.1} Definition}  We say that a type $p \in
S(A)$ is \ub{almost realized} in $B$ (or $B$ almost realizes $p$) if
for all $\varepsilon > 0$ there exists $b_\varepsilon \in
B,b_\varepsilon \models p^{[\varepsilon]}$.
\enddefinition
\bigskip

\proclaim{\stag{MOM.4B} Claim}  [$({\frak C},\bold d)$ compact]  Let $M \in
K^c_1$ be such that every 1-type over a subset of $M$ of cardinality
less than $\lambda$ is almost realized in $M$ ($\lambda$ infinite).
Then $M$ is $(D,\lambda)$-homogeneous.
\endproclaim
\bigskip

\demo{Proof}  Let $p \in S^1_D(A),A \subseteq M,|A| < \lambda$.
Choose by induction $b_n \models p^{[\frac{1}{2^n}]},b_n \in M,\bold
d(b_{n+1},b_n) \le \frac{1}{2^{n-1}}$.

Why is this possible?  Let $q_n = p(x) 
\cup\{\bold d(x,b_n) \le \frac{1}{2^n}\}$.
As $b_n \models p^{[\frac{1}{2^n}]}$, some $b^* \in {\frak C}$
realizes $q_n$ (see \scite{MOM.4}) and let $q^*_n = \text{ tp}_D(b^*/A
\cup \{b_n\})$, so $q^*_n \in S^1_D(A \cup \{b_n\})$, and by the
assumption there exists $b_{n+1} \in M$ realizing
$(q^*_n)^{[\frac{1}{2^{n+1}}]}$.  So $b_{n+1} \models
p^{[\frac{1}{2^{n+1}}]},\bold d(b_{n+1},b_n) \le \frac{1}{2^n} +
\frac{1}{2^{n+1}} \le \frac{1}{2^{n-1}}$, as required.

Now, $\langle b_n:n < \omega\rangle$ is a Cauchy sequence, let $b \in
M$ be its limit ($M \in K^c_1$, so is complete).  Now $b \models
p^{[\varepsilon]}$ for all $\varepsilon > 0$, so $b \models p$ (see
\scite{MOM.3}(1)), and we are done. 
\enddemo
\bigskip

\proclaim{\stag{MOM.4C} Corollary}  Let $M \in K^c_1$ be
non-$(D,\lambda)$-homogeneous.  Then there exists $A \subseteq M,|A| <
\lambda,p \in S^1_D(A)$ and $\varepsilon > 0$ such that
$p^{[\varepsilon]}$ is omitted in $M$.
\endproclaim
\bigskip

\demo{Proof}  This is just restating \scite{MOM.4B}; but we prefer
this form for later use.
\enddemo
\bigskip

\remark{\stag{MOM.4D} Remark}  We will not use the notion of pseudo
homogeneity (Definition \scite{MOM.4A}), in this paper (as a
postoriori all the models will turn out to be
$(D,\lambda)$-homogeneous), but it is interesting to point out what 
non-categorical classes can look like.  See \scite{4n.23} later. 
\endremark
\bigskip

\proclaim{\stag{MON.2s.1} Claim}   If $A \subseteq {\frak C},p \in 
\bold S_D(\text{\rm mcl}(A))$ and $c \in {\frak C}$ realizes $p
\restriction A$ \ub{then} it realizes $p$.
\endproclaim
\bigskip

\demo{Proof}  If $\varphi(x,a_1,\dotsc,a_m) \in p$ ($\varphi \in
\Delta$ of course, parameters not suppressed) then for each $n <
\omega$ we can choose $b^n_\ell \in A$ (for $\ell =
1,\dotsc,m)$ such that $\bold d(a_\ell,b^n_\ell) < \frac{1}{n+2}$.
Hence $\varphi(x,b^n_1,\dotsc,b^n_m)^{<1/n+2>}$ belongs to $p$ as
$\varphi(x,a_1,\dotsc,a_m)$ implies it, hence it belongs to $p \restriction A$.

So if $c \in {\frak C}$ realizes $p \restriction A$ then $\models
\varphi(c,b^n_1,\dotsc,b^n_m)^{<1/n+2>}$ for each $n$.  If ${\frak C} \models
\neg \varphi[c,a_1,\dotsc,a_m]$ then some $\psi,\psi =
\psi(y,x_1,\dotsc,x_m) \in \Delta$ and $\varphi(y,x_1,\dotsc,x_n)$ are
contradictory and ${\frak C} \models \psi[c,a_1,\dotsc,a_n]$.  So for
all $n$, the tuples $cb^n_1,\dotsc,b^n_m$ lies in both 
$(\varphi^{<\frac{1}{n+2}>})^{\frak C}$
and $(\psi^{<\frac{1}{n+2}>})^{\frak C}$, therefore 
$\varphi^{<\varepsilon>} \cap
\psi^{<\varepsilon>} \ne \emptyset$ for all $\varepsilon > 0$, which
contradicts $\varphi$ and $\psi$ being a contradictory pair, see
\scite{MOM.3A}(4). 
\enddemo
\bigskip

\head {\S4 Stability in momspaces} \endhead  \resetall \sectno=4
 \spuriousreset
\bigskip

We define a topological version of stability.  The intuition behind
the definition is that there may be many types, but the density of the
space of types is small.  Our definition generalizes
Iovino's stability for Banach spaces, see \cite{Io99} for Hausdorff
cats and continuous theories, it coincides with definitions given in
\cite{BY05}, \cite{BeUs0y} respectively.  Note that for
elementary homogeneous class (i.e. discrete metric), this definition
coincides with the usual one (so certainly for an elementary class,
i.e., ${\frak C}$ saturated).  For a non-discrete metric the classical
$\lambda$-stability of $D$ (counting types, as in \cite{Sh:3}) is
stronger than the topological relative we define here, and is equivalent
if and only if $\lambda = \lambda^{\aleph_0}$.  Stability in our sense
(i.e., $\lambda$-stable for some $\lambda$) is equivalent to stability
for $D$, but for a specific $\lambda$ (e.g. $\lambda = \aleph_0$) the
notions differ.
\bigskip

\demo{\stag{STAB.0} Hypothesis}  $({\frak C},\bold d)$ is a momspace for
$\Delta = \Delta_{\frak C}$.
\enddemo
\bigskip

\definition{\stag{STAB.1} Definition}  1) A momspace $({\frak C},\bold
d)$ is called \ub{$0^+-\lambda$-stable} if for all $A \subseteq {\frak
C}$ of cardinality $\lambda$ there exists $B \subseteq {\frak C}$ of
cardinality $\lambda$ such that each $p \in \bold S_D(A)$ is realized
in mcl$(B) = \bar B$ (i.e., the topological closure of $B$).
\enddefinition
\bigskip

\proclaim{\stag{STAB.2} Claim}  [$({\frak C},\bold d)$ compact]  
For a compact momspace ${\frak C}$, the following are equivalent:
\mr
\item "{$(A)$}"   ${\frak C}$ is $0^+-\lambda$-stable
\sn
\item "{$(B)$}"  for each $A \subseteq {\frak C},|A| \le \lambda$, there exists
$B,|B| \le \lambda$ such that for each $\varepsilon >0$ and $p \in
\bold S_D(A),p^{[\varepsilon]}$ is realized in $B$ ($B$ almost realizes
all types over $A$, see Definition \scite{MOM.4A.1})
\sn
\item "{$(C)$}"  there are no 
$A,|A| \le \lambda,\langle p_i:i < \lambda^+\rangle$ a
sequence of members of $\bold S_D(A)$, i.e.,  complete 
$D$-types over $A$ and $\varepsilon > 0$ such that
$p^{[\varepsilon]}_i \cap p^{[\varepsilon]}_j = \emptyset$ for all
$i,j < \lambda^+$.
\endroster
\endproclaim
\bigskip

\demo{Proof}  $(A) \Rightarrow (C)$.  

Assume (A).  Suppose (C) fails,
so we have $\langle p_i:i < \lambda^+\rangle$ over $A$ as there.  By
(A) we can find $B$ of cardinality $\le \lambda$ 
such that each $p_i$ is realized in mcl$(B) = \bar B$, by,
say, $a_i$.  As $\{a_i:i < \lambda^+\}$ form an $\varepsilon$-net (by
the nature of the $p_i$'s), it is obvious that density of $\bar B$ is
at least $\lambda^+$, contradiction ($|B|=\lambda$).
\mn
$(C) \Rightarrow (B)$.

Let $A$ be given.  We construct $B_n$ by induction:

$-B_0 = A$

$-B_{n+1}$ realizes all complete $D$-types over $B_n$ up to $\frac{1}{n+1}$,
i.e.
\nl

\hskip20pt if $p \in S_D(B_n)$ then for some $a \in B_{n+1}$,
tp$(a,B_n)^{[\frac{1}{n+1}]}({\frak C}) \cap p^{[\frac{1}{n+1}]}
({\frak C}) \ne \emptyset$

$-|B_n| = \lambda$

$B = \cup B_n$ is obviously as required in (B).

How is $B_{n+1}$ constructed?  Let $\varepsilon = 1/2(n+1)$ and let
$\langle p^k_i:i < \lambda
\rangle$ be a maximal $\varepsilon$-disjoint set of complete types of
$\kappa$-tuples over $B_n$, (I.e. $(p^k_i)^{[\varepsilon]} \cap
(p^k_j)^{[\varepsilon]} = \emptyset$ for all $i,j$)

$B_{n+1} = B_n \cup \{\bar a^k_i:\bar a^k_i$ realizes $p^k_i$ for some
$i < \lambda,\kappa < \omega\}$.

Now each complete $k$-type $p$ over $B_n$ is $\frac{1}{n+1}$-realized in
$B_{n+1}$, as there is $i < \lambda$ such that
$(p^k_i)^{[\varepsilon]}({\frak C}) \cap p^{[\varepsilon]}({\frak C})$
is non-empty, let $\bar b$ be in the intersection.  We can replace
$\bar b$ by any $\bar b'$ realizing tp$(b,B_n)$ hence \wilog \, $\bold
d(\bar b,a^k_i) < \varepsilon$.  Also we can find $\bar c$ realizing
$p$ such that $\bold d(\bar b,\bar c) < \varepsilon$.  So $\bar a^k_i$
witnesses $p$ is $\frac{1}{n+1}$-realized in $B_{n+1}$.
\mn
(B) $\Rightarrow$ (A).

Let $A$ be given.  Define $B_0 = A,B_{n+1}$ almost realizes all types
over $B_n,B_n \subseteq B_{n+1},|B_n| = \lambda$.  Let $B = \dbcu_{n < \omega}
B_n$.  Let $p \in \bold S_D(A)$.  Pick by induction on $n \ge 1,a_n
\in B_n$ such that $a_n \models p^{[\frac 1n]}$ and 
$\bold d(a_n,a_{n-1}) \le \frac {1}{2^n}$ if $n>1$ 
(possible by \scite{MOM.4} as $B_n$
almost realizes all types over $\dbcu_{i < n} B_n$).  $\langle
a_n:n <\omega \rangle$ is obviously a Cauchy sequence, let $a$ be its
limit in $\bar B$ and we are done by \scite{MOM.3}(1) 
as $a \models p^{[\frac 1n]}$ for all $n$.
\enddemo
\bigskip

\demo{\stag{STAB.2A} Observation}  1) If ${\frak C}$ is
$\lambda$-stable (in the sense of $D$, see \cite{Sh:3}), then it 
is $0^+-\lambda$-stable.
\nl
2) If ${\frak C}$ is $0^+-\lambda$-stable and $\lambda =
\lambda^{\aleph_0}$ then ${\frak C}$ is $\lambda$-stable.
\enddemo
\bigskip

\demo{Proof}  1) Trivial.
\nl
2) Recall $|\bar B| = |B|^{\aleph_0}$.
\enddemo
\bn
Recall
\definition{\stag{STAB.2s.2} Definition}  A (partial) type $q$ \ub{splits}
over a set $A$ if there are $\varphi(\bar x,\bar a),\psi(\bar x,\bar
b) \in q,\varphi(\bar x,\bar y),\psi(\bar x,\bar y)$ contradictory,
and tp$(\bar a/A) = \text{\rm tp}(\bar b/A)$.
\enddefinition
\bigskip

\demo{\stag{STAB.3} Observation}  For a complete type $p,p^{[\varepsilon]}$
splits over $A \Leftrightarrow$ there exist $\varphi(\bar x,\bar
a),\psi(\bar x,\bar b) \in p,\bold d_1(\varphi(\bar x,\bar
y),\psi(\bar x,\bar y)) > 2 \varepsilon$, 
tp$(\bar a/A) = \text{\rm tp}(\bar b/A)$.
\enddemo
\bigskip

\demo{Proof}  $p^{[\varepsilon]}$ splits over $A \Leftrightarrow$
there exist $\varphi^{[\varepsilon]}(\bar x,\bar
a),\psi^{[\varepsilon]}(\bar x,\bar b) \in p^{[\varepsilon]}$ (where
$\varphi(\bar x,\bar a),\psi(\bar x,\bar b) \in p)$ such that $\bar a
\equiv_A \bar b$ and $(\varphi^{[\varepsilon]}(\bar x,\bar
y),\psi^{[\varepsilon]}(\bar x,\bar y))$ is a contradictory pair [so
$(\varphi,\psi)$ is a $2\varepsilon$-contradictory pair].
\enddemo
\bigskip

\demo{\stag{STAB.3s.1} Observation}  1) The type $p$ splits over $A$
\ub{iff} for some $\varepsilon >0$, the type $p^{[\varepsilon]}$
splits over $A$.
\enddemo
\bigskip

\demo{Proof}  $p$ splits over $A$ iff some contradictory pair
$(\varphi,\psi)$ witnesses this, now pick $\varepsilon = \frac{\bold
d_1(\varphi,\psi)}{3}$ and use \scite{STAB.3}.
\enddemo
\bigskip

\demo{\stag{STAB.3.4} Notation}  We write $\bar a \equiv_A \bar b$ for
tp$(\bar a/A) = \text{\rm tp}(\bar b/A)$.
\enddemo
\bigskip

\proclaim{\stag{STAB.4} Lemma}  Let ${\frak C}$ be $0^+-\aleph_0$-stable.  
\ub{Then} for any 
$B \subseteq {\frak C},p \in \bold S_D(B)$ and $\varepsilon > 0,
p^{[\varepsilon]}$ does not split over a finite subset of $B$.
\endproclaim
\bigskip

\demo{Proof}  Suppose $p^{[\varepsilon]}$ splits over every finite subset
of its domain.   We construct finite sets $A_n$ for $n < \omega$ and
elementary maps $F_\eta$ for $\eta \in {}^{\omega >}2$ as follows:
\mr
\item  $A_0 = \emptyset$
\sn
\item   $A_{n+1} = A_n \cup \{\bar a_n,\bar b_n\}$, where there
are $\varphi_n(\bar x,\bar y),\psi_n(\bar x,\bar y)$
such that $\varphi_n(\bar x,\bar a_n),\psi_n(\bar x,\bar b_n)$ exemplify
$p^{[\varepsilon]}$ splits over $A_n$, i.e. $\bar a_n \equiv_{A_n}
\bar b_n,\varphi_n(\bar x,\bar a_n),
\psi_n(\bar x,\bar b_n) \in p^{[\varepsilon]}$ and
$\varphi^{[\varepsilon]}_n({\frak C},\bar a_n) \cap
\psi^{[\varepsilon]}_n({\frak C},\bar b_n) = \emptyset$
\sn
\item   for $\eta \in {}^n 2,F_\eta:A_n \rightarrow {\frak C}$ is an
elementary mapping
\sn
\item  for $\eta \in {}^n 2,F_{\eta {}^\frown \,\langle 0\rangle}
(\bar a_n) = F_{\eta {}^\frown \,\langle 1 \rangle} (\bar b_n)$ and 
$F_{\eta {}^\frown \,\langle 0 \rangle},F_{\eta {}^\frown\,\langle 1 \rangle}$ 
extend $F_\eta$.
\ermn
The construction is straightforward.  Now denote for $\eta
\in {}^\omega 2,F_\eta = \dbcu_{n < \omega} F_{\eta \upharpoonright n},
p^*_\eta = F_\eta(p),A = \cup\{\text{\rm Rang}(F_\eta):\eta
\in {}^{\omega >} 2\}$ (so it is countable) and choose \footnote{this uses
``$D$ is good"} $p_\eta,p^*_\eta \subseteq p_\eta \in \bold S_D(A)$.  
Obviously, $\eta \ne \nu \in
{}^\omega 2 \Rightarrow p_\eta^{[\varepsilon]} \cap
p^{[\varepsilon]}_\nu = \emptyset$, contradicting
$0^+-\aleph_0$-stability by \scite{STAB.2} by an implication not using
compactness.
\enddemo
\bigskip

\proclaim{\stag{STAB.2s.0} Claim}  Assume 
$({\frak C},\bold d)$ is $0^+-\aleph_0$-stable.  \ub{Then}
\mr
\item "{$(a)$}"  if $p \in \bold S_D(B)$ then $p$ does not split over some
countable $A \subseteq B$
\ermn
\mr
\item "{$(b)$}"  $D$ is stable (in the sense of \cite{Sh:3}), $\kappa(D)
\le \aleph_1$.
\endroster
\endproclaim
\bigskip

\remark{Remark}  Note that it is close but not as in first order;
there may be $\aleph_0$ exceptions.
\endremark
\bigskip

\demo{Proof}  
\mn
\ub{Clause (a)}:

By the previous claim for every $\varepsilon >0$ for some finite
$B_\varepsilon \subseteq A,p$ does not
$\varepsilon$-split over $B_\varepsilon$.  Let $B = \cup\{B_{1/(n+1)}:n
< \omega\}$, so $B$ is a countable subset of $A$ and by
\scite{STAB.3s.1} and the obvious monotonicity of
non-$\varepsilon$-splitting, $p$ does not split over $B$. 
\mn

\ub{Clause (b)}:  Follows from (a).
\bigskip

\definition{\stag{ST.1} Definition}  Given an uncountable indiscernible
set $\bold I \subseteq {}^m {\frak C}$ and a set $A \subseteq {\frak
C}$, define the average type of $\bold I$ over $A$, Av$(A,\bold I)$ as
follows:

$$
\align
\text{Av}(A,\bold I) = \{\varphi(\bar x,\bar a):&\varphi(\bar x,\bar y)
\in \Delta,\ell g(\bar x) = m,\bar a \in A, \text{ and for infinitely many} \\
 &\bar c \in \bold I,\varphi(\bar c,\bar a) \text{ holds}\}.
\endalign
$$
\enddefinition
\bigskip

\demo{\stag{ST.2} Fact}  If ${\frak C}$ is stable, then any
indiscernible sequence is an indiscernible set.
\enddemo
\bigskip

\demo{Proof}  Standard.

We often say ``$\bold I$ is indiscernible" meaning an indiscernible
sequence, which is the same as indiscernible set.
\enddemo
\bigskip

\proclaim{\stag{ST.3} Claim}  Let ${\frak C}$ be $0^+-\aleph_0$-stable.
\nl
1) If $\bold I \subseteq {}^m {\frak C}$ is indiscernible uncountable,
$A \subseteq {\frak C}$ a set, then {\rm Av}$(A,\bold I) \in
S^m_D(A)$.
\nl
2) If $A,\bold I$ are as in (1), then

$$
\align
\text{Av}(A,\bold I) = \{\varphi(\bar x,\bar a):&\varphi(\bar x,\bar y)
\in \Delta,\ell g(\bar x) = m,\bar a \in M \text{ and for all but} \\
  &\text{ countably many } \bar c \in \bold I \text{ does } 
\varphi(\bar c,\bar a) \text{ hold}\}.
\endalign
$$
\endproclaim
\bigskip

\demo{Proof}  Using the standard argument,
one shows that for a given contradictory pair $(\varphi,\psi) =
(\varphi(\bar x,\bar z),\psi(\bar x,\bar z))$ and $\bar a \in M$, one
of the sets $\{\bar c \in \bold I:\models \varphi(\bar c,\bar
a)\},\{\bar c \in \bold I:\models \psi(\bar c,\bar a)\}$ is finite
(otherwise, let $\varepsilon = \bold d_1(\varphi,\psi)$, and construct
$2^{\aleph_0} \, \varepsilon$-distant types over a countable set,
contradictory $0^+-\aleph_0$-stability; in fact, one constructs
$2^{\aleph_0}$ pairwise distinct $(\varphi,\psi)$-types, i.e., types
mentioning only $\varphi$ and $\psi$).

Now given a formula $\varphi(\bar x,\bar a)$ over $M$, if $\{\bar c
\in \bold I:\varphi(\bar c,\bar a)$ holds$\}$ is infinite, then for
each $\psi(\bar x,\bar z)$ such that $(\varphi,\psi)$ is
contradictory, the set $J_\psi = \{\bar c \in \bold I:\psi(\bar c,\bar
a)$ holds$\}$ is necessarily finite, so taking the union of $J_\psi$
over all such $\psi$, we obtain a countable set of exceptions

$$
J = \cup J_\psi = \{\bar c \in \bold I:\models \neg \varphi(\bar
c,\bar a)\}.
$$
\mn
This completes the proof of clause (2).  For clause (1), let
$\varphi(\bar x,\bar a)$ be a formula over $A,\ell g(\bar x) = m$.  If
for uncountably many $\bar c \in \bold I,\varphi(\bar c,\bar a)$ holds,
then $\varphi(\bar x,\bar a) \in \text{ Av}(A,\bold I)$, otherwise for
some $\psi(\bar x,\bar a)$ such that $(\varphi,\psi)$ is a
contradictory pair, $\psi(\bar c,\bar a)$ holds for uncountably many
$\bar c \in \bold I$, so $\psi(\bar x,\bar a) \in \text{ Av}(A,\bold
I)$.  Clearly, only one of the two options above is possible, so (1) follows.
\enddemo
\bn
\ub{Discussion}:  Why countable and not finite in the definition of
averages?  Even if the majority satisfies
$\varphi(\bar x,\bar a)$, for each $\varepsilon$ there can be finitely many
$\bar c \in \bold I$ such that 
${\frak C} \models \psi(\bar c,\bar a)$ and $(\varphi,\psi)$ are
$\varepsilon$-contradictory, and this finite number can increase when
$\varepsilon$ goes to 0. 
\enddemo
\bigskip

\proclaim{\stag{ST.4} Lemma}  If $M$ is 
$(D,\aleph_1)$-homogeneous, $p \in \bold S^m(M)$
\ub{then} for some uncountable $\bold I \subseteq {}^m M$, we have $p
= \text{\rm Av}(M,\bold I)$.
\endproclaim
\bigskip

\demo{Proof}  Let $B$ be as in \scite{STAB.2s.0},
clause (a).  Let $m=1$ for simplicity.  Choose $a_\alpha \in M$ 
realizing $p \restriction B \cup\{a_\beta:\beta < \alpha\}$ by
induction on $\alpha < \omega_1$.
Now $\bold I = \langle a_\alpha:\alpha < \omega_1\rangle$ is
indiscernible over $B$ by \cite[I,\S2]{Sh:c}.  

If $q = \text{Av}(M,\bold I) \ne p$ still $q \in \bold S_D(M)$ (by 
\scite{ST.3}(1)) and
we can find $\varphi(x,\bar b) \in q,\psi(x,\bar b) \in p$ and they
are contradictory.  So $u = \{\alpha < \omega_1:{\frak C} \models
\varphi[a_\alpha,\bar b]\}$ is infinite let $v \subseteq u$ be of
cardinality $\aleph_0$, let $v \subseteq \alpha(*) < \omega_1$ and
choose $a'_\beta \in M(\beta \in [\alpha(*),\omega_1)$ realizing $p
\restriction (B \cup \{a_i:i < \alpha(*)\} \cup \bar b \cup
\{a'_\gamma:\gamma \in (\alpha(*),\beta)\}$.
\nl
Easy contradiction: for a given contradictory pair $(\varphi,\psi)$
all but finitely many elements of $\bold I$ have to ``make a choice", see
also \scite{ST.3}(1).
\enddemo
\bigskip

\definition{\stag{STAB.2s.1} Definition}  1) A momspace 
$({\frak C},\bold d)$ is called $(\mu,*)$-superstable if 
given $\langle M_i:i < \omega \rangle$ an
increasing chain of $(D,\mu)$-homogeneous models $(M_i \in K^c_1$, 
of course), mcl$(\dbcu_{i < \omega} M_i)$ is $(D,\mu)$-homogeneous.
\nl
2) We omit $\mu$ if this holds for every $\mu$ large enough.
\nl
3) Fully $*$-superstable means for every $\mu > |\tau_{\frak C}|
   + \aleph_0$.
\enddefinition
\bigskip

\remark{Remark}  This definition generalizes superstability for
${\frak C}$ a saturated model of a first order theory, and 
${\frak C}$ a homogeneous monster.
\endremark
\bn
The following claim will be mainly of interest for us when $\mu = \aleph_0$:
\proclaim{\stag{STAB.5} Claim}  Let $({\frak C},\bold d)$ be
$0^+-\aleph_0$-stable and compact.  Then $({\frak C},\bold d)$ is 
$(\mu^+,*)$-superstable for every $\mu \ge \aleph_0$.
\endproclaim
\bigskip

\demo{Proof}  Let $\langle M_n:n < \omega \rangle$ be an increasing
sequence of $(D,\mu^+)$-homogeneous models, and assume $p \in
\bold S(A),A \subseteq \overline{\dbcu_{n < \omega} M_n}$ of cardinality
$\le \mu$ is not realized in $M_\omega = \overline{\bigcup M_n}$.  

By increasing $A$ and the $0^+-\aleph_0$-stability 
(i.e., for every $\varepsilon > 0$ trying to build a tree 
$\langle p_\eta:\eta \in {}^{\omega >} 2\rangle$ of 
$\varepsilon$-contradictory types) \wilog \, $p$ has a unique
extension in $\bold S_D(M_\omega)$, call it $q$.
By \scite{MOM.4C} as $({\frak C},\bold d)$ is compact we can add
that for some
$\varepsilon >0,p^{[\varepsilon]}$ is not realized in $M_\omega$.
Without loss of generality $A \subseteq  \text{\rm mcl}
(A \cap \dbcu_{n < \omega} M_n)$ hence by \scite{MON.2s.1} \wilog \,
$A \subseteq \dbcu_{n < \omega} M_n$ (as
$p$ is determined by its restriction to $A \cap (\dbcu_{n < \omega} M_n))$.

Now by \scite{STAB.4} there is a finite $B
\subseteq A$ over which $q \restriction \cup\{M_n:n < \omega\}$ 
does not $(\varepsilon/5)$-split.  Let $n
< \omega$ be such that $B \subseteq M_n$ \wilog \, 
$q \restriction M_n$ does not split over $A \cap M_n$ (by increasing
$A$ and \scite{STAB.2s.0}(a)). 
Let $q_n =q \restriction M_n$ and let $A_n =
M_n \cap A$.  As $M_n$ is $(D,\mu^+)$-homogeneous, by Lemma
\scite{ST.4}, there is an
uncountable indiscernible sequence $\bold I_n$ in $M_n$ with
Av$(M_n,\bold I_n) = q_n$; \wilog \,
$|\bold I_n| = \mu^+$, and $\bold I_n$ is indiscernible over $A$ (not just
$A_n$!) (as for each finite type $\bar b$ from $A$, tp$(\bar b/\bold
I_n)$ does not split over a countable subset, so we can remove a
subset of $\bold I_n$ of cardinality $\mu$).  Now
as elements of $\bold I_n$ do not realize $p^{[\varepsilon]}$, for
some formula $\varphi(\bar x,\bar a) \in p$ (really $\bar x$ is a singleton)
there exists $\psi(\bar x,\bar z)$ such that $\bold d_1(\varphi,\psi)
> \varepsilon$ and $\forall \bar c \in \bold I_n,\psi(\bar c,\bar a)$.

For $k < \omega$ let 
$\bar c_k \in \bold I_n$ be pairwise distinct.  We can find $\bar a' 
\in M_n$ which realizes tp$(\bar a,A_n \cup B \cup \cup\{\bar c_k:k
<\omega\})$.  So by clause $(c)$ of \scite{STAB.2s.0} we know that for
all but countably many $\bar c \in \bold I_n$ we have ${\frak C}
\models \psi[\bar c,\bar a']$ (as this happens for $\{\bar c_k:k <
\omega\} \subseteq I_n$), hence $\psi(\bar x,\bar a') \in
q_n \subseteq q$.
\enddemo
\bn
Now we obtain: for some $m> n,\bar a \in M_m$ hence
$(\varphi(\bar x,\bar a),\psi(\bar x,\bar a'))$ witness that
$q \restriction \cup\{M_\ell:\ell < \omega\}$ does $\varepsilon$-split
over $B$, hence $q$ does, which is a contradiction to the choice of
$B$.    \hfill$\square$
\bn
We shall now proceed to proving an analogue of density of isolated
types.
As the ``right" notion of a type in our context seems to
be an $\varepsilon$-neighborhood of a complete type, the assumption of
``non-isolated"  will not be enough for us.

\definition{\stag{B.10.1} Definition}  We say that $M \in K_1$ 
$(<\varepsilon)$-omits $p(\bar x)$, a type over $A \subseteq M$, when for
no $\zeta \in [0,\varepsilon)_{\Bbb R}$ and $\bar b \in M$, does
$\bar b$ realize $p^{[\zeta]}$.
\nl
\enddefinition
\bigskip

\definition{\stag{4n.9} Definition}  1) We say that a formula
$\psi(\bar x,\bar b)$ pseudo $(\varepsilon,\zeta)$ isolates a type $p(\bar
x)$ if $\psi^{<\zeta>}(\bar x,\bar b) \models p^{[\varepsilon]}(\bar x)$
(note that the roles of $\varepsilon,\zeta$ are not symmetric and the
different notions of approximation!).  In other words, if
${\frak C} \models \psi^{<\zeta>}[\bar a,\bar b]$ then for
some $\bar a' \in p({\frak C})$ we have $\bold d(\bar a',\bar a)
\le \varepsilon$ (so $\bar a'$ realizes $p,\bar a$ realizes
$\psi^{<\zeta>}(\bar x,\bar b)$).
\nl
2) We say that $A$ is a pseudo $(< \varepsilon)$-support for 
$p(\bar x)$ or $A$ pseudo $(< \varepsilon)$-supports $p$ 
when there is a consistent
$\psi(\bar x,\bar b),\bar b \subseteq A$ and positive $\zeta_1,\zeta_2$ such
that $\psi^{<\zeta_1>}(\bar x,\bar b) \models
p^{[\varepsilon-\zeta_2]}(\bar x)$.  So $\psi(\bar x,\bar b)$ psuedo 
$(\varepsilon-\zeta_2,\zeta_1)$-isolates $p(\bar x)$. 
\nl
3) We say that $A$ really $(< \varepsilon)$-omits $p(\bar x)$ if it 
does not pseudo $(< \varepsilon)$-support $p(\bar x)$.
\enddefinition
\bigskip

\proclaim{\stag{4n.9.1} Claim}  1) If $M \in K^c_1 \, (< \varepsilon)$-omits
$p(\bar x),p(\bar x) \in \bold S^m(A),A \subseteq M$ 
\ub{then} $M$ really $(< \varepsilon)$-omits $p(\bar x)$.
\nl
2) If $p(\bar x)$ is a type over $M$ and $M$ really 
$(< \varepsilon)$-omits $p(\bar x)$, then $M(< \varepsilon)$-omits $p$.
\nl
\endproclaim
\bigskip

\demo{Proof}  1) Assume $M$ is a pseudo $(< \varepsilon)$-support for $p(\bar
x)$, i.e., there exist $\psi(\bar x,\bar b)$ over $M$ and
$\zeta_1,\zeta_2 > 0$ such that $\psi^{<\zeta_1>}(\bar x,\bar b)
\models p^{[\varepsilon - \zeta_2]}(\bar x)$.  ${\frak C} \models
\exists x \psi(\bar x,\bar b),M \in K^c_1$, so for some $\bar a \in
M,{\frak C} \models \psi^{<\zeta_1>}(\bar a,\bar b)$ (see
\scite{MOM.suB}), therefore
$p^{[\varepsilon-\zeta_2]}(\bar a)$ holds, and we have
$p^{(<\varepsilon)}(\bar x)$ is realized in $M$.
\nl
2) Easier: assume $\bar a \models p^{[\varepsilon-\zeta]}(\bar x)$ for
$\zeta > 0,\bar a \in M$.  Then the formula ``$\bar x = \bar a$" is
over $M$ and pseudo $(\varepsilon - \zeta_2,\zeta_1)$-isolates $p(\bar
x)$ for $\zeta_2 = \zeta_1 = \frac{\zeta}{3}$.
\enddemo
\bigskip

\definition{\stag{4n.13} Definition}  1)  We say 
$\varphi(\bar x,\bar b)$ strictly
$(\varepsilon,\zeta)$-isolates a type $p$ if $\varphi(\bar x,\bar b) \in p$ and
$\varphi(\bar x,\bar b)$ pseudo $(\varepsilon,\zeta)$-isolates $p$.
\nl
2) We say $\varphi(\bar x,\bar
b)$ is strictly $(\varepsilon,\zeta)$-isolating over $A$ when
\mr
\item "{$(a)$}"  $\bar b \subseteq A$
\sn
\item "{$(b)$}"  if $\varphi(\bar x,\bar b) \in p \in \bold S^{\ell g(\bar
x)}_D(A)$ then $\varphi(\bar x,\bar b)$ strictly isolates $p$.
\ermn
3) We say that $p \in S^m_D(A)$ is strictly
$(\varepsilon,\zeta)$-isolated if some $\varphi$ strictly
$(\varepsilon,\zeta)$-isolates it.
\nl
3A) ``Strictly $\varepsilon$-isolate" means ``for some $\zeta >
0$, strictly $(\varepsilon,\zeta)$-isolate".
\nl
4) We say that $p \in \bold S^m_D(A)$ is strictly isolated when for every
$\varepsilon > 0$ for some $\varphi(\bar x,\bar a) \in p$ and some
$\zeta > 0$ the formula $\varphi(\bar x,\bar a)$ 
strictly $(\varepsilon,\zeta)$-isolates the type $p$ (i.e., $p$ is
$\varepsilon$-strictly isolated for all $\varepsilon > 0$).
\enddefinition
\bigskip

\proclaim{\stag{4n.14} Claim}  [$({\frak C},\bold d)$ is compact]

Assume that
\mr
\item "{$(a)$}"  $p \in \bold S^m_D(A)$ or just $p$ is an $m$-type closed
under conjunctions
\sn
\item "{$(b)$}"  $\psi(\bar x,\bar b) \in p$
\sn
\item "{$(c)$}"  $\varepsilon > 0$ and $\zeta \ge 0$.
\ermn
\ub{Then} one of the following occurs
\mr
\item "{$(\alpha)$}"  there is a pair $(\psi_1(\bar x,\bar
y),\psi_2(\bar x,\bar y))$ of formulas and a sequence $\bar b^*$ from $A$ such
that $\bar b \triangleleft \bar b^*,
\ell g(\bar b^*) = \ell g(\bar y)$ such that
$\psi(\bar x,\bar b) \wedge \psi_1(x,\bar b^*),\psi^{<\zeta>}
(\bar x,\bar b) \wedge \psi_2(\bar x,\bar b^*)$ are 
$\varepsilon$-contradictory (and both consistent, of course), 
hence for no $\bar a,{\frak C}
\models (\exists \bar x)(\bold d(\bar x,\bar a) \le \varepsilon/2
\wedge \psi(\bar x,\bar b) \wedge \psi_1(\bar x,\bar b^*))$ and
${\frak C} \models (\exists \bar x)[\bold d(\bar x,\bar a) \le
\varepsilon/2 \wedge \psi^{<\zeta>}(\bar x,\bar b) \wedge \psi_2(\bar
x,\bar b^*)]$
\sn
\item "{$(\beta)$}"  for every $\bar a'$ such that ${\frak C}
\models \psi^{<\zeta>}[\bar a',\bar b]$ there is a sequence $\bar a''$
realizing $p$ such that $\bold d(\bar a',\bar a'') \le \varepsilon$
(so $\psi(x,\bar b)$ strictly $(\varepsilon,\zeta)$-isolates $p$, see
\scite{4n.13}).
\endroster
\endproclaim
\bigskip

\demo{Proof}  We can assume that clause $(\beta)$ fails and let $\bar a'$ 
exemplify it.  So for every $\bar a''$ realizing $p$ we have
$\bold d(\bar a'',\bar a') > \varepsilon$.

Let $q(\bar y) = \text{\rm tp}(\bar a',A)$ so $\psi^{<\zeta>}(\bar
y,\bar b) \in q$ as ${\frak C} \models \psi^{<\zeta>}(\bar a',\bar b)$,  
and let $r(\bar x,\bar y) =
p(\bar x) \cup q(\bar y) \cup \{\bold d(\bar x,\bar y) \le
\varepsilon\}$.  If $r(\bar x,\bar y)$ is consistent, so is $r(\bar
x,\bar a')$, and any $\bar a''$ realizing $r(\bar x,\bar a')$ is as
required in clause $(\beta)$, contradicting our assumption.
So $r(\bar x,\bar y)$ is inconsistent. 
As $({\frak C},\bold d)$ is compact and $p(\bar x),q(\bar y)$ are
closed under conjunctions,
$\psi(\bar x,\bar b) \in p(\bar x),\psi^{<\zeta>}(\bar y,\bar b) \in q(\bar y)$
and as we can add dummy
variants, there are $\bar b^* \subseteq A,\bar b \triangleleft \bar
b^*$ and $\psi_1(\bar x,\bar b^*) \in p(\bar x),\psi_2(\bar
y,\bar b^*) \in q(\bar y)$ such that $\{\psi(\bar x,\bar b) \wedge
\psi_1(\bar x,\bar b^*),\psi^{<\zeta>}(\bar x,\bar b) \wedge
\psi_2(\bar y,\bar b^*),\bold d(\bar x,\bar y) \le
\varepsilon\}$ is contradictory.  \hfill$\square$

So we get clause $(\alpha)$.
\enddemo
\bigskip

\proclaim{\stag{4n.12} An isolation Claim}  [$({\frak C},\bold d)$ is
$0^+-\aleph_0$-stable and compact]
\nl
1) If $\bar a \subseteq A$ and $\varphi(\bar x,\bar a)$ is consistent
and $\varepsilon > 0$ \ub{then} we can find $\varphi_1(\bar x,\bar
a_1)$ and $\zeta >0$ such that $\bar a_1 \subseteq A$ and $\varphi(\bar
x,\bar a) \wedge \varphi_1(\bar x,\bar a_1)$ is strictly
$(\varepsilon,\zeta)$-isolating over $A$, see Definition \scite{4n.13}.
\nl
1A) Similarly omitting $\zeta$ getting strictly $\varepsilon$-isolating.
\nl
2) The set of strictly isolated $p \in S^m_D(A)$ is dense, i.e. for
 every $\varphi(\bar x,\bar a)$ with $\bar a \in A$, there exists $p
\in S^m_D(A),\varphi(\bar x,\bar a) \in p,p$ is strictly isolated.
\endproclaim
\bigskip

\demo{Proof of \scite{4n.12}}  Part (1A) is restating Part (1).  
Also part (2) follows from part (1) by choosing $\varphi_n(\bar x,\bar
a_n)$ such that $\varphi(\bar x,\bar a) \wedge \varphi_1(\bar x,\bar
a_1) \wedge \ldots \wedge \varphi_n(\bar x,\bar a_n)$ is
$\frac{1}{n+1}$-isolating over $A$ (iterating 1A) and applying
compactness of ${\frak C}$.  So we concentrate on proving part (1) .

We can choose $\zeta_n > 0$ for $n < \omega$ such that,
e.g. $\Sigma\{\zeta_n:n <  \omega\} \le \varepsilon/5$ 
\nl

Assume that $\psi(\bar x,\bar b),\bar b \subseteq A$ is a
counterexample.  Now we choose $\langle \psi_\eta(\bar x,\bar
a_\eta):\eta \in {}^n 2 \rangle$ by induction on $n$ such that
\mr
\item "{$\boxtimes$}"  $(a) \quad \bar a_\eta \subseteq A$
\sn
\item "{${{}}$}"  $(b) \quad \psi_\eta(\bar x,\bar a_\eta)$ is
consistent
\sn
\item "{${{}}$}"  $(c) \quad \psi_{<>}(\bar x,\bar a_{<>}) = \psi(\bar
x,\bar b)$
\sn
\item "{${{}}$}"  $(d) \quad$ if $\nu {}^\frown \langle 0 \rangle,\nu
{}^\frown \langle 1 \rangle \in {}^n 2$ then $\psi_{\nu {}^\frown
\langle 0 \rangle}(\bar x,\bar a_{\nu {}^\frown \langle 0 \rangle}),
\psi_{\nu {}^\frown \langle 1 \rangle}(\bar x,\bar a_{\nu
{}^\frown \langle 1 \rangle})$ are 
\nl

\hskip20pt $\varepsilon$-contradictory
\sn
\item "{${{}}$}"  $(e) \quad$ if $\eta = \nu {}^\frown \langle 0
\rangle \in {}^n 2$ then $\psi_\eta(\bar x,\bar a_\eta) \models
\psi_\nu(\bar x,\bar a_\nu)$
\sn
\item "{${{}}$}"  $(f) \quad$ if $\eta = \nu {}^\frown \langle 1
\rangle \in {}^n 2$ then $\psi_\eta(\bar x,\bar a_\eta) \models
\psi^{<\zeta_n>}_\nu(\bar x,\bar a_\nu)$.
\ermn
By \scite{4n.14} there is no problem to carry the definition, i.e.,
having $\psi_\nu(\bar x,\bar a_\nu)$, clause $(\beta)$ of
\scite{4n.14} cannot hold (with $\psi_\nu,\bar a_\nu$ here standing
for $\psi,\bar b$ there) as ``$\psi(\bar x,\bar a)$ is a
counterexample".  Hence clause $(\alpha)$ there holds, let us
choose $\psi_{\nu {}^\frown <\ell>}(\bar x,\bar a_{\nu {}^\frown
<\ell>})$ for $\ell=0,1$.  

Now let $\xi_n = \Sigma\{\zeta_m:m \in [n,\omega)\}$, so clearly
\mr 
\item "{$(*)$}"  if $n(1) < n(2) < \omega$ and
$\eta_\ell \in {}^{n(\ell)}2$ for $\ell=1,2$ and $\eta_1
\triangleleft \eta_2$ then $\psi^{<\xi_{n(2)}>}_{\eta_2}(\bar
x,\bar a_{\eta_2}) \models \psi^{<\xi_{n(1)}>}_{\eta_1}(\bar x,
\bar a_{\eta_1})$
\nl
[Why?  By \scite{4n.10} using clauses (e) + (f) of $\boxtimes$ we get
$\psi_{\eta_2} \models \psi^{<\xi_{n(1),n(2)}>}_{\eta_1}$, where
$\xi_{n(1),n(2)} = \Sigma\{\zeta_m:m \in [n(1),n(2))\}$.  Now use
\scite{4n.10} again.]
\ermn
Now let $C = \cup\{\bar a_\eta:\eta \in {}^{\omega >}2\}$.
\nl
Hence
\mr
\item "{$(*)$}"  for $\eta \in {}^\omega 2$ the set
$\{\psi^{<\xi_n>}_{\eta \restriction n}(\bar x,\bar a_{\eta
\restriction n}):n < \omega\}$ is consistent, hence is included in some
$p_\eta \in \bold S(C)$
\sn
\item "{$(*)$}"  if $\nu {}^\frown \langle \ell \rangle \triangleleft
\eta_\ell \in {}^\omega 2$ for $\ell=0,1$ then 
$p^{[\varepsilon/5]}_{\eta_1}(\bar x) \cup
p^{[\varepsilon/5]}_{\eta_2}(\bar x)$ is inconsistent.
\nl
[Why?  By $\boxtimes(d)$ and the choice of $\zeta_n$], a contradiction
to $0^+-\aleph_0$-stability.  
\nl
${{}}$   \hfill$\square_{\scite{4n.12}}$
\endroster
\enddemo
\bn
\centerline {$* \qquad * \qquad *$}
\bn
The following is not used at present but clarifies non-categoricity.
Recall Definition \scite{MOM.4A} and Claim \scite{MOM.4B}.  Note that
\scite{MOM.4C} says that a non-$(D,\lambda)$-homogeneuos model omits
some $p^{[\varepsilon]}$.  Here we clarify what happens in the case of
pseudo $(D,\lambda)$-homogeneous non-$(D,\lambda)$-homogeneous model.
\proclaim{\stag{4n.23} Claim}  Assume $M \in K^c_1,\lambda > \aleph_0 +
|\tau_{\frak C}|$ and $M$ is not $(D,\lambda)$-homogeneous.  Then $(*)$ or
 $(**)$
\mr
\item "{$(*)$}"  $(a) \quad N \prec^1_\Delta M,|N| < \lambda$
\sn
\item "{${{}}$}"  $(b) \quad p \in \bold S_D(N)$
\sn
\item "{${{}}$}"  $(c) \quad$ if ${\frak C}$ is $0^+-\mu$-stable for
some $\mu \in [\aleph_0 + |\tau_{\frak C}|,\lambda)$ then $p$ has a 
\nl

\hskip25pt unique extension in $\bold S_D(M)$
\sn
\item "{${{}}$}"  $(d) \quad \varepsilon > 0$ and $M$ omits
 $p^{<\varepsilon>}$
\sn
\item "{$(**)$}"  $(a) \quad$ for every $B \subseteq M,|B| < \lambda$
and $q \in \bold S_D(B)$ and $\zeta > 0$, the type $q^{<\zeta>}$
\nl

\hskip25pt  is realized in $M$ (so $M$ is pseudo $\lambda$-saturated)
\sn
\item "{${{}}$}"  $(b) \quad N \prec^1_\Delta M,|\tau_{\frak C}| +
\aleph_0 \le \|N\| < \lambda$
\sn
\item "{${{}}$}"  $(c) \quad p \in \bold S_D(N)$
\sn
\item "{${{}}$}"  $(d) \quad$ if ${\frak C}$ is $0^+-\mu$-stable for
some $\mu \in [|\tau_{\frak C}| + \aleph_0,\lambda)$ then $p$ has a
unique
\nl

\hskip25pt  extension in $\bold S_D(M)$
\sn
\item "{${{}}$}"  $(e) \quad \varepsilon > 0,p^{[\varepsilon]}$ is
 omitted by $M$
\sn
\item "{${{}}$}"  $(f) \quad$  for every $n$, for some
$c_n,\varepsilon_n,\zeta_n$ we have $1/(n+1) > \varepsilon_n >
\zeta_n$
\nl

\hskip25pt $> 0,c_n \in M$ realizes $p^{<\varepsilon_n>}$ and $\bold
d(c_n,p^{<\zeta_n>}({\frak C})) \ge 10 \times \varepsilon_n - \zeta_n$.
\endroster
\endproclaim
\bigskip

\demo{Proof}   Clearly there is $A \subseteq M$ such that 
$|A| < \lambda$ and $p \in \bold S^1_D(A)$ is omitted.
Let $\mu = |A| + |\tau_{\frak C}| + \aleph_0$ so $\mu < \lambda$.  If
${\frak C}$ is $0^+-\mu'$-stable for some $\mu' \in [|\tau_{\frak C}|
+ \aleph_0,\lambda)$, easily \wilog \, $p$ has a unique extension in
$\bold S_D(M)$ and $A = |N|,N \prec^1_\Delta M$.  If for some
$\varepsilon > 0,p^{<\varepsilon>}$ is also omitted by $M$,  the case $(*)$
holds.  So we may assume $(*)$ fails, in other words, clause (a) of
$(**)$ holds.

Let $n^* < \omega,\varepsilon^* = \frac 12$ (any $\varepsilon^* > 0$
works).  We are going to find $c_{n^*},\varepsilon_{n^*},\zeta_{n^*}$
as required in clause (f) above.  First, we try to choose $b_n
\in M$ by induction on $n < \omega$ such that
\mr
\item "{$\circledast_n$}"  $(a) \quad b_n \in M$
\sn
\item "{${{}}$}"  $(b) \quad b_n$ realizes $p^{<\varepsilon^*/(n+1)>}$
\sn
\item "{${{}}$}"  $(c) \quad$ if $n=m+1$ then $\bold d(b_n,b_m) \le 10
\times \varepsilon^*/2^n$.
\ermn
For $n=0$ there is $b_0 \in M$ realizing $p^{<\varepsilon^*/1>} =
p^{<\varepsilon^*>}$ by clause (a) of $(**)$.  So we can begin.

\bn
\ub{Point 1}:   We cannot succeed to choose $\langle b_n:n < \omega
\rangle$.
\nl
Why?  Suppose we have succeeded.  Then
$\langle b_n:n < \omega \rangle$ is a Cauchy sequence and therefore
converges to some $b^* \in M$.  We will show that $b^* \models p$.
If $\varphi(x,\bar a) \in p$, then for each $n,{\frak C} \models
\varphi^{<\varepsilon^*/(n+1)>}(b_n,\bar a)$ so there is $\bar a_n
\in {}^{\omega >}{\frak C},b'_n \in {\frak C}$ such that
\mr
\item "{$\boxtimes$}"   $(a) \quad \bold d(b'_n,b_n) \le
\varepsilon^*/(n+1)$
\sn
\item "{${{}}$}"  $(b) \quad \bold d(\bar a_n,\bar a) \le
\varepsilon^*/(n+1)$
\sn
\item "{${{}}$}"  $(c) \quad \models \varphi[b'_n,\bar a_n]$.
\ermn
Now 

$$
\langle b'_n:n < \omega \rangle \text{ converges to } b^*
$$

$$
\langle \bar a_n:n < \omega \rangle \text{ converges to } \bar a
$$
\mn
hence $\models \varphi[b^*,\bar a]$.  So $b^* \models p,b* \in M$, a 
contradiction.
\enddemo
\bn
\ub{Point 2}:

So we are stuck in some $n=m+1$ so let $c_{n^*} := b_m,
\varepsilon_{n^*} = \varepsilon^*/2^n,\zeta_{n^*} =
\varepsilon^*/2^{n+2}$.  If the demand in (f) of $(**)$ fails, then
there is $b'_m \in {\frak C}$ realizing $p^{<\zeta_{n^*}>}(\bar x) \cup
\{\bold d(x,b_m) \le 10 \times \varepsilon_{n^*} - \zeta_{n^*}\}$ hence
$p^{<\zeta_{n^*}>}(x) \cup \{\bold d(x,b_m) \le 10 \times
\varepsilon_{n^*} - \zeta_{n^*}\}$
is consistent hence it is contained in some $q_n \in \bold S(N \cup\{b_m\})$.

So for every $\zeta >0$ (we use $\zeta$ small enough) there is $b_n
\in M$ realizing $q^{<\zeta>}(x)$ (recall we are assuming
$(**)(a)!$).  So $b_n$ realizes $(p^{<\zeta_{n^*}>})^{<\zeta>}$ 
hence $p^{<\zeta_{n^*}+ \zeta>}$ hence if
$\zeta$ is small enough, $p^{<\varepsilon^*/(n+1)>}$.  Also $\bold
d(b_n,b_m) \le (10 \times \varepsilon_{n^*} - \zeta_{n^*}) + \zeta + \zeta
< 10 \times \varepsilon_{n^*} = 10 \times \varepsilon^*/2^n$ 
(because if $(a',b')$ realized $\bold d^{<\zeta>}(x,y) \le \xi$ 
then $\bold d(a',b') \le \xi + \zeta +
\zeta$).

So $b_n$ is as required in $\circledast_n (a)-(c)$ above, so we could
have continued choosing the $b_n$.  \hfill$\square_{\scite{4n.12}}$
\bigskip

\head {\S5 Ehrenfeucht-Mostowski models} \endhead  \resetall \sectno=5
 \spuriousreset
\bigskip

In this section we adapt the technique of constructing
Ehrenfeucht-Mostowski models to our context.  The reader should have a
look at chapter 7 of \cite{Sh:c} for the basic definitions ($\Phi$
proper, etc.).  The basic idea is the following: we start with ${\frak
C}$ in vocabulary $\tau$.  Adding skolem functions, we obtain
vocabulary $\tau'$.  Choosing an indiscernible sequence and taking its
type (its EM - ``blueprint") $\Phi$, for each order type $J$ we can
construct EM$(J,\Phi)$ (like in chapter 7 of \cite{Sh:c}), which will
be an elementary submodel of ${\frak C}$ expanded to $\tau'$, therefore
its restriction to $\tau$, EM$_\tau(J,\Phi)$ is an elementary submodel of
${\frak C}$, although not necessarily complete.  Taking the
completion, we obtain a model in $K^c_1$. 
Adding more structure to the language we can make it
$(D,\lambda)$-homogeneous, and more, see below.  

Given a vocabulary $\tau^*$ with skolem functions, and a
$\tau^*$-diagram of indiscernibles $\Phi$ (EM-blueprint), we denote
for each order-type $I$, the EM-model (the $\tau^*$-skolem hull of a
sequence $\langle a_i:i \in I \rangle$) by EM$_{\tau^*}(I,\Phi)$ or
EM$(I,\Phi)$ if $\tau^*$ is clear from the context.  We denote by
EM$_{\tau_0}(I,\Phi)$ the restriction of EM$(I,\Phi)$ to the vocabulary
$\tau_0 \subseteq \tau^*$.

Let ${\frak C}$ be a momspace. Let $\tau$ be the vocabulary of ${\frak C}$.  
It is clear that for any $\tau^*$ (with skolem functions) expanding
$\tau$, a $\tau^*$-diagram of indiscernibles $\Phi$ (in ${\frak C}$ expanded
to $\tau^*$), $I$ an order, we can think of
EM$_\tau(I,\Phi)$ as an elementary submodel of ${\frak C}$, so 
EM$_\tau(I,\Phi) \prec_{\Bbb L(\tau({\frak C}))} {\frak C}$.  This is
not necessarily true for the completion, but 
$\overline{\text{EM}_\tau(I,\Phi)} \prec^1_\Delta {\frak C}$ by
\scite{MOM.2.19}. 
\bigskip

\proclaim{\stag{EM.1} Claim}  Let $({\frak C},\bold d)$
be a momspace, $\tau$ the
vocabulary of ${\frak C},|\tau| \le \aleph_0,\tau' \supseteq
\tau,\tau'$ with Skolem functions, $\Phi'$ a $\tau'$-blueprint.
\nl
0A) For every linear order $J$, {\rm EM}$_\tau(J,\Phi') \in K_1$
\nl
and {\rm mcl}$(\text{\rm EM}_\tau(J,\Phi')) \in K^c_1$.  

0B) If $J_1 \subseteq J_2$ then {\rm EM}$_\tau(J_1,\Phi')
\prec^1_\Delta \text{\rm EM}_\tau(J_2,\Phi')$; moreover {\rm
EM}$_\tau(J_1,\Phi') \prec \text{\rm EM}_\tau(J_1,\Phi')$ and
{\rm mcl}$(\text{\rm EM}_\tau(J_1,\Phi')) \prec^1_\Delta \text{\rm
mcl}(\text{\rm EM}_\tau(J_2,\Phi'))$.
\nl

1) There exists $\tau^*$ expanding $\tau',|\tau^*| = 2^{\aleph_0}$ and
a $\tau^*$-diagram $\Phi^*$ such that for each finite order $J$,
EM$_\tau(J,\Phi^*)$ is $(D,\aleph_1)$-homogeneous.
\nl
2) If ${\frak C}$ is $(\aleph_1,*)$-superstable, then $\Phi^*$ as in (1) works
for all orders $J$, but we have to take the closure, 
i.e., mcl(EM$_\tau(J,\Phi^*))$ which $\in K^c_1$ is 
$(D,\aleph_1)$-homogeneous for all $J$.
\nl
3) If ${\frak C}$ is $0^+-\omega$-stable, then $\tau^*$ in (1) and (2)
can be
chosen of cardinality $\aleph_1$.
\endproclaim
\bigskip

\demo{Proof}   0) (A),(B) straight.
\nl
1)  Choose for $i < \aleph_1,\tau_i,\Phi_i,|\tau_i| 
= 2^{\aleph_0}$ with skolem functions expanding $\tau'$ 
increasing continuous such that for each
$\tau_i$-type $p$ over a finite subset of the skeleton of
EM$(I,\Phi_i)$, say $a_1,\dotsc,a_n$, there exists a function symbol
$f_p$ in $\tau_{i+1}$ such that $f_p(a_1,\dotsc,a_n)$ realizes $p$.

More precisely, we do the following:  for any consistent set $p$ of
formulas of the form $\varphi = \varphi(x,y_1,\dotsc,y_n) \in \tau_i$
($y_1,\ldots,y_n$ are the parameters; some of the $y_i$'s may be 
dummy variables) such that $p$ is closed under conjunctions and for every
$\varphi \in p,\exists x \varphi(x,y_1,\dotsc,y_n) \in \Phi_i$, we add
 a function symbol $f_p$ to $\tau_{i+1}$ such that for every $\varphi \in
p$ the following formula is in $\Phi_{i+1}:\exists x
\varphi(x,y_1,\dotsc,y_n) \rightarrow
\varphi(f_p(y_1,\dotsc,y_n),y_1,\dotsc,y_n)$.
\nl
Now let $\Phi^* = \dbcu_{i < \aleph_1} \Phi_i$ and let $M = \text{\rm
EM}_\tau(J,\Phi^*)$ for some finite $J$.  Choose $A \subseteq M$
countable, $p \in S(A)$.  As $A$ is countable, it can be viewed as a
countable subset of EM$(J,\Phi_i)$ for some $i,p$ is a type over the finite
skeleton, so realized in EM$(J,\Phi_{i+1})$, therefore in $M$, as
required.
\nl
2) By induction on $|J|$.   
We just need to show that for an increasing sequence of linear orders
$J_i,\dbcu_i \text{\rm EM}(J_i,\Phi) = 
\text{\rm EM}(\dbcu_i J_i,\Phi)$ and this is clear as
elements of EM$(J,\Phi)$ have finite character, i.e., use only
finitely many elements of the skeleton $J$.  Of course, we then have
to take metric closure.
\nl
3) Let $J_n$ be a linear order with $n$ elements.
Similarly to (1), we choose by induction on $i < \omega_1$
countable $\tau_i \subseteq \tau^*$ increasing continuous, closed
under skolem functions, and $\Phi_i$ such that each type over
EM$(J_n,\Phi_i)$ is almost realized in EM$(J_n,\Phi_{i+1})$: we choose a
countable set $B$ which almost realizes all types over EM$(J_n,\Phi_i)$,
and for each such $p$ and for each $\kappa$ we have 
$f_{p,\kappa} \in \Phi_{i+1}$ such
that $f_{p,\kappa}(\bar a)$ is $\frac 1k$-close to a realization of
$\varphi$ for each $\varphi \in p$.
\enddemo
\bigskip

\proclaim{\stag{EM.2} Corollary}  If ${\frak C}$ is
$0^+-\omega$-stable, it has a $(D,\aleph_1)$-homogeneous model
in all uncountable density characters.
\endproclaim
\bigskip

\demo{Proof}  Let $\lambda > \aleph_0$.
Consider mcl(EM$(\lambda,\Phi^*)),\Phi^*$ as in
\scite{EM.1}(3) (note that $|\Phi^*| \le \aleph_1$) and use
\scite{EM.1}(2) + \scite{STAB.5}.
\enddemo
\bn
\margintag{EM.3}\ub{\stag{EM.3} Discussion}:  If $({\frak C},\bold d)$ is
$0^+-\aleph_0$-stable, does it have a $(D,\lambda)$ homogeneous model in
every $\lambda$?  By \scite{UNI.1} this follows from categoricity,
which is good enough for our purposes.
\bigskip

\head {\S6 Embeddings, isomorphisms and categoricity} \endhead  \resetall \sectno=6
 \spuriousreset
\bigskip

In this section we introduce notions of $\varepsilon$-embedding and
$\varepsilon$-isomorphism which are weaker than isometry.  
This will lead us to the notion of weak uncountable categoricity that
we investigate in \S8.
\bigskip

\demo{\stag{ISO.0S} Convention}  Models are from $K=K_1$.
\enddemo
\bigskip

\definition{\stag{ISO.1} Definition}  For two metric structures in the
same vocabulary $\tau$ and $\varepsilon \ge 0$ we say
\mr
\item  $f:M_1 \rightarrow M_2$ is an $\varepsilon$-embedding if for every
$\Delta$-formula $\varphi,\bar a \in M_1,M_1 \models \varphi(\bar a)
\Rightarrow M_2 \models \varphi^{[\varepsilon]}(\bar a)$
\sn
\item  $f:M_1 \rightarrow M_2$ is an $\varepsilon$-isomorphism if it is
an $\varepsilon$-embedding which is one-to-one and onto
\sn
\item  $M_1,M_2$ are $\varepsilon^+$-isomorphic if there exists a
$\zeta$-isomorphism $f_\zeta:M_1 \rightarrow M_2$ for all $\zeta > 
\varepsilon$
\sn
\endroster
\enddefinition
\bigskip

\demo{\stag{ISO.2} Observation}  1)  $0$-embedding is a regular notion of
(isometric) embedding, $0$-isomorphisms is regular isomorphism (in
particular isometry).
\nl

2) If there exists a $\zeta$-isomorphism $f_\zeta:M_1 \rightarrow
M_2$ for all $\zeta > \varepsilon$, then there exists a
$\zeta$-isomorphism $g_\zeta:M_2 \rightarrow M_1$ for all $\zeta >
\varepsilon$ (so clause (3) of the definition above makes sense).
\enddemo
\bigskip

\demo{Proof}  Clear.
\enddemo
\bigskip

The following definition is the central one.
\definition{\stag{ISO.3} Definition}  Let $({\frak C},\bold d)$ be a
momspace, $\varepsilon \ge 0,\lambda$ a cardinal.
\nl
1) We say ${\frak C}$ is $\varepsilon^+$-categorical in $\lambda$ if
every two complete $M_1,M_2 \in K^c_1$ of density $\lambda$ are
$\varepsilon^+$-isomorphic.
\nl
2) We say ${\frak C}$ is categorical in $\lambda$ if every two complete
$M_1,M_2 \in K^c_1$ of density $\lambda$ are isomorphic.
\nl
3) We say that ${\frak C}$ is possibly categorical
($\varepsilon^+$-categorical) if it is categorical 
($\varepsilon^+$-categorical) in some $\lambda > \aleph_0$.
\nl
4) We say that ${\frak C}$ is weakly uncountably categorical
(wu-categorical) if the following holds: for each $\varepsilon >0$
there exists a cardinal $\lambda$ such that ${\frak C}$ is
$\varepsilon^+$-categorical in $\lambda$.
\enddefinition
\bigskip

\demo{\stag{ISO.4} Observation}  1) If $\varepsilon \ge \zeta$ then
$\zeta^+$-categoricity implies $\varepsilon^+$-categoricity (for a
specific $\lambda$).
\nl
2) Possible $0^+$-categoricity implies weak uncountable
categoricity.
\nl
3) Categoricity implies all the other notions (for a specific $\lambda$).
\enddemo
\bigskip

\proclaim{\stag{OS} Theorem}  [$({\frak C},\bold d)$ compact]
\nl
Let $K = K^c_1({\frak C})$ be a wu-categorical
momspace with countable language.  Then ${\frak C}$ is $0^+-\aleph_0$-stable.
\endproclaim
\bigskip

\demo{Proof}  Let $\tau$ be $\tau(K),\tau'$ is $\tau$ expanded with skolem
functions, $\Phi$-proper for $K$.  So $\tau'$ is countable.
\enddemo
\bigskip

\demo{\stag{OS.1} Subclaim}  Under these assumptions, let $I$ be a
well-ordered set, $M_0 = \overline{\text{EM}_\tau(I,\Phi)}$, $A
\subseteq M_0$ is countable $\varepsilon > 0$.  Then 
each $\varepsilon$-disjoint set ${\Cal P}$ of types from $\bold S_D(A),
\frac{\varepsilon}{2}$-realized in $M_0$ (so
$p_1,p_2 \in {\Cal P} \Rightarrow p^{[\varepsilon]}_1 \cup
p^{[\varepsilon]}_2$ is contradictory and for each $p \in 
{\Cal P},p^{[\frac{\varepsilon}{2}]}$ is realized in $M_0$) is countable.
\enddemo
\bigskip

\demo{Proof of the Subclaim}  If not, let $\langle p_i:i < \omega_1
\rangle$ be $\varepsilon$-disjoint types over
$A,p^{[\frac{\varepsilon}{2}]}_i$ realized in $M_0$ by $\bar b_i$.
Pick $\bar b^0_i \in$ {\rm EM}$(I,\Phi),\bold d(\bar b^0_i,\bar b_i) \le
\frac{\varepsilon}{100}$.  Without loss of generality $A = \text{\rm
EM}(J,\Phi)$ for $J \subseteq I,|J| \le \aleph_0$.  As $J$ is well
ordered, by the standard argument, there are uncountably many
$b^0_i$'s satisfying the same type over $A$, but $b^0_i \models
p^{[\varepsilon]}_i$ and $p^{[\varepsilon]}_i,p^{[\varepsilon]}_j$ are
contradictory for $i \ne j$, a contradiction.  \hfill$\square_{\scite{OS.1}}$
\enddemo
\bn
Now we prove the theorem.  Assuming ${\frak C}$ is not
$0^+-\omega$-stable, we get $A \subseteq {\frak C}$ countable and
$\langle p_i:i < \omega_1 \rangle \, \varepsilon$-disjoint types over
$A$ for some $\varepsilon > 0$ (remember \scite{STAB.2}).  Let $\lambda$
be such that ${\frak C}$ is $\delta^+$-categorical in $\lambda$ for
$\delta << \varepsilon$.  Now apply the usual argument: choose $M_1
\in K$ of density $\lambda$ which includes $A$ and $\langle b_i:i <
\omega_1\rangle$ realizations of $\langle p_i:i < \omega_1 \rangle$,
and on the other hand consider $M_0 =
\overline{\text{EM}(\lambda,\Phi)}$.  Applying $f:M_1 \rightarrow 
M_0$ which is a $\delta_1$-embedding $\delta_1 <
\frac{\varepsilon}{2}$, we get that $\langle f(b_i):i
< \omega_1\rangle$ contradict \scite{OS.1}.  \hfill$\square_{\scite{OS}}$
\bigskip

\proclaim{\stag{MOM.EXIST} Corollary}  Let ${\frak C}$ and $K$ be as in
\scite{OS}, then $K$ has a $(D,\aleph_1)$-homogeneous model
in $\lambda$ (recall this means of density $\lambda$) 
for all $\lambda > \aleph_0$.
\endproclaim
\bigskip

\demo{Proof}  By \scite{OS} and \scite{EM.2}.
\enddemo
\bigskip

\head {\S7 Uni-dimensionality} \endhead  \resetall \sectno=7
 \spuriousreset
\bn

The following notion was explored in \cite{Sh:3} but not defined there:
\definition{\stag{UD.1} Definition}  A (good) finite diagram 
$D$ is \ub{uni-dimensional} if for some regular $\lambda$ there is no
$(D,\lambda)$-homogeneous model of $K$ which is not
$\lambda^+$-homogeneous in cardinality $\ge \lambda^+$.
\enddefinition
\bn
In \cite{Sh:3} it is essentially proven (see \cite[\S6]{Sh:3}) that:
\proclaim{\stag{UD.2} Theorem}  Assume $D$ is stable.  Then the
following are equivalent:
\mr
\item  $D$ is not uni-dimensional
\sn
\item  there is some regular $\lambda$ such that there are maximally
$(D,\lambda)$-homogeneous models of arbitrary large cardinalities 
\sn
\item  for all large enough regular $\lambda < \mu$, there is a
$(D,\lambda)$-homogeneous model $M$ and $\langle a_i:i < \mu
\rangle,\langle b_i:i < \lambda \rangle$ mutually indiscernible
sequences in $M$ such that $\langle b_i:i < \lambda \rangle$ is a
maximal indiscernible sequence in $M$, i.e., can not be extended in
$M$
\sn
\item  there is a cardinal $\lambda$ and a model $M$ of cardinality
$\lambda$ which is $(D,\aleph_1)$-homogeneous, but not
$(D,\lambda)$-homogeneous. 
\endroster 
\endproclaim
\bigskip

\remark{\stag{UD.2A} Remark}  In our context we will say ``${\frak C}$
is uni-dimensional" or ``$K$ is uni-dimensional", meaning that $D$ is.
\endremark
\bigskip

\demo{\stag{UD.3} Reminder}  1) For an indiscernible set $\bold I \subseteq 
{\frak C}$ of cardinality $> |\tau_{\frak C}| + \aleph_0$ 
and a set $A \subseteq {\frak C}$, we define

$$
\align
\text{Av}(\bold I,A) = \{\varphi(\bar x,\bar a):&\bar a \in A, 
\text{ infinitely many elements} \\
  &\text{ of } \bold I \text{ satisfy } \varphi(\bar x,\bar a)\}.
\endalign
$$
\mn
We call this set the \ub{average type} of $\bold I$ over $A$.  (See
\scite{STAB.2s.0} ``all but finitely many" is wrong, see
\cite{Sh:863} on ``cutting indiscernibles", the weaker version fits here.)
\nl
2) For stable ${\frak C}$, for any indiscernible $\bold I,|\bold I| 
> |\tau_{\frak C}| + \aleph_0$ and set $A$, Av$(\bold I,A)$ is a
complete type, see \scite{ST.3}(1).
\nl
3) Let $\bold I = \langle \bar a_i:i < \alpha \rangle$, where $\langle \bar
a_i:i \le \alpha\rangle$ is indiscernible.  Then $\bar a_\alpha
\models \text{\rm Av}(\bold I,\cup \bold I)$.
\nl
4) Let $\bold I$ be indiscernible, $\bold I = 
\langle \bar a_i:i < \alpha \rangle$
and let $\bar a_\alpha \models \text{\rm Av}(\bold I,\cup \bold I)$.  
Then $\langle \bar a_i:i \le \alpha \rangle$ is indiscernible.
\nl
5) It follows from (3) + (4) that $\bold I = \langle \bar a_i:i < \alpha
\rangle \subseteq M$ is a maximal indiscernible sequence (set) in $M$
iff Av$(\bold I,\cup \bold I)$ is omitted in $M$.
\nl
6) $\varphi(x,\bar a_{i_1},\dotsc,\bar a_{i_n}) \in \text{\rm
Av}(\bold I,\cup \bold I)$
for $\bold I = \langle \bar a_i:i < \delta\rangle$ ($\delta$-limit ordinal) iff
$\bar a_j \models \varphi(x,\bar a_{i_1},\dotsc,\bar a_{i_n})$ 
for some/every $j \notin \{i_1,\dotsc,i_n\}$.
\enddemo
\bigskip

\proclaim{\stag{UNI} Theorem}  1) Let $({\frak C},\bold d)$ be a
momspace, $\tau({\frak C})$ countable, $0^+$-categorical in
$\lambda,\lambda > \aleph_0$.  Then ${\frak C}$ is uni-dimensional.
\nl
2) The same is true if ${\frak C}$ is wu-categorical.
\endproclaim
\bigskip

\demo{Proof}  1) If not, choose $0 < \theta_1 << \theta_2$ and let $M$
be a model $M \in K,|M| = \theta_2,M$ is $\theta_1$-homogeneous,
$\langle a_i:i < \theta_2 \rangle,\langle b_i:i < \theta_1 \rangle$
mutually indiscernible, $\langle a_i:i < \theta_1 \rangle$ cannot be
extended in $M$, so (denoting $\bold I = \langle a_i:i < \theta_2
\rangle,\bold J = \langle b_i:i < \theta_1 \rangle$) 
Av$(\bold J,\cup \bold J)$ is omitted in $M$. 

We now expand the language by a predicate $P$ for $J$ and skolem
functions, call the new vocabulary $\tau'$.  Let $T' = \text{\rm
Th}_{\tau'}(M')$ (where $M'$ is $M$ in the expanded language).

Note that
\mr
\item "{$\circledast_0$}"  $T' \models ``P$ is a $\tau$-indiscernible
set",i.e., for every $\tau$-formula, $T'$ implies that any two tuples
from $P$ behave the same.
\ermn
The type $p = \text{\rm Av}_\tau(\bold J, \cup \bold J)$ is omitted 
in $M$, therefore \wilog \, by \scite{MOM.4C},
for some $\varepsilon,p^{[\varepsilon]}$ is omitted.  If we choose
$\theta_1,\theta_2$ carefully enough $(\theta_2 >> \theta_1)$ in the
beginning, then by the Erd\"os-Rado theorem (as in the proof of
\cite[VIII,5.3]{Sh:c}, or using \cite{BY03a}(1.2)) we can choose a diagram
of indiscernibles (EM-blueprint) $\Phi$ in vocabulary $\tau'$ such
that for any $\mu$, denoting the skeleton of $M'_0 :=
\,\text{\rm EM}(\mu,\Phi)$ by $I' = \langle a'_i:i < \mu \rangle$, we have
\mr
\item "{$\circledast_1$}"  $\bold I'$ is a $\tau'$-indiscernible sequence
(set), moreover, it is $\tau'$-indiscernible over $P^{M'_0}$
\sn
\item "{$\circledast_2$}"  for each $n < \omega$, for some
$i_1,\dotsc,i_n < \theta_2,a'_0,\dotsc,a'_{n-1}$ has the same
$\tau'$-type as $a_{i_1},\dotsc,a_{i_n}$.
\ermn
Now:
\mr
\item "{$\circledast_3$}"  let $M_0 = \text{ EM}_\tau(\mu,\Phi)$, then
$M_0 \prec {\frak C}$, so mcl$(M_0) \in K = K^c_1$ (as on the one hand
$\tau'$ has skolem functions, and on the other hand $M'_0$ does not
realize $\tau$-types over $\emptyset$ that were not realized in $M$,
so $M_0$ is a $D$-model)
\sn
\item "{$\circledast_4$}"  $M'_0 \models T'$ (skolem functions, so
$M'_0 \equiv M'$)
\sn
\item "{$\circledast_5$}"  $P^{M'_0}$ is a $\tau$-indiscernible set
(by $\circledast_4$ and $\circledast_0$ above)
\sn
\item "{$\circledast_6$}"  $P^{M'_0}$ is countable.  Why?  Each $b \in
P^{M'_0}$ is of the form $\sigma(a'_{i_1},\dotsc,a'_{i_n})$ for some
$\tau'$-term $\sigma$.  But as $\bold I'$ is 
$\tau'$-indiscernible over $P^{M'_0}$,
each such $b$ depends only on $\sigma$, and there are countably many
$\tau'$-terms ($\tau$ is countable, and so is $\tau'$).
\ermn
Denote $\bold J' = P^{M'_0}$, a $\tau$-indiscernible set.
Denote $p' = \text{\rm Av}_\tau(\bold J', \cup \bold J')$.  Then 
\nl
$[p']^{[\varepsilon]}$ is omitted in $M_0$.
\nl
Why?  Pick $\sigma(a'_{i_j},\dotsc,a'_{i_n}) \in M'_0$.  Let
$j_1,\dotsc,j_n$ be such that $a'_{i_1},\dotsc,a'_{i_n} \equiv
a_{j_1},\dotsc,a_{j_n},\sigma(a_{j_1},\dotsc,a_{j_n}) \in M'$ does not
realize $p^{[\varepsilon]}$, so for some $\varphi(x) \in p,M' \models
\bold d_1(\sigma(\bar a_j),\varphi) \ge \varepsilon$.
\nl
Note: $\varphi(x) = \varphi(x,\bar c),\bar c \in P^{M'}$.  Call
that $\varphi(x,\bar c) \in p \Leftrightarrow M' \models
\varphi(d,\bar c)$ for some/all $d \in P^{M'},d \cap \bar c =
\emptyset$ (as $P^{M'} = J$ is an indiscernible set, see
\scite{UD.3}(6)).  So $M' \models
``\exists \bar c \in P$ such that $[\forall d \in P \backslash \bar
c,\varphi(d,\bar c)] \and [\bold d_1(\sigma(\bar a_j),\varphi(x,\bar c)) \ge
\varepsilon]"$.  Therefore, $M'_0$ satisfies the same formula with
$\sigma(\bar a'_i)$, which obviously means that $\sigma(\bar a'_i)$ does not
satisfy $[p']^{[\varepsilon]}$, as required.

We have finished now: let $\mu = \lambda$, so in $\bar M_0 =
\overline{\text{EM}_\tau(\lambda,\Phi)}$ we have a countable
indiscernible set whose average is omitted (as $(p')^{[\varepsilon]}$
is omitted in $M_0$), so $\bar M_0$ is a
non-$(D,\aleph_1)$-homogeneous model in ${\frak K}$ of density
$\lambda$, but by \scite{MOM.EXIST} we have a
$(D,\aleph_1)$-homogenous model of density $\lambda > \aleph_0$.  So
categoricity in $\lambda$ fails, moreover, $0^+$-categoricity fails,
as $\bar M_0$ is at least $\frac{\varepsilon}{2}$-distant from any
$(D,\aleph_1)$-homogeneous model. 
\nl
2) Repeat the proof of (1), and at the end choose $\lambda$ in which,
say, $(\frac{\varepsilon}{2})$-categoricity holds, and get the same
contradiction. 
\enddemo
\bigskip

\proclaim{\stag{UNI.1} Claim}  Let ${\frak C}$ be a wu-categorical
momspace, $K = K^c_1({\frak C})$.
\nl
1) There exists a $(D,\lambda)$-homogeneous model in $K$ for all
$\lambda > \aleph_0$. 
\nl
2) Each $(D,\aleph_1)$-homogeneous model in $K$ of density $\lambda$
is $(D,\lambda)$-homogeneous.
\nl
3) If $K$ is $0^+$-categorical in $\lambda > \aleph_0$, then each 
$K$-model of density $\ge \lambda$ is $(D,\lambda)$-homogeneous.
\endproclaim
\bigskip

\demo{Proof}  1) By \scite{MOM.EXIST} and uni-dimensionality. 
\nl
2) By the equivalence \scite{UD.2}.
\nl
3) Otherwise by a L\"owenheim-Skolem argument we will get a
 non-$(D,\lambda)$-homogeneous model of density $\lambda$, and
 together with (1) this will lead to a contradiction.
\enddemo
\bigskip

\head {\S8 The main theorem} \endhead  \resetall \sectno=8
 \spuriousreset
\bigskip

\demo{\stag{CATEG.0} Hypothesis}  ${\frak C}$ is a compact momspace 
with countable vocabulary.
\enddemo
\bigskip

\proclaim{\stag{CATEG} Theorem}  Assume $K = K^c_1({\frak C})$ wu-categorical.  
\ub{Then} $K$ is $\lambda$-categorical for
all $\lambda > \aleph_0$, moreover, any model of $K$ of
density $\lambda > \aleph_0$ is $(D,\lambda)$-homogeneous.
\endproclaim
\bigskip

\demo{Proof}  Suppose not, so ${\frak C}$ is $0^+-\aleph_0$-stable, 
$(\aleph_1,*)$-superstable,  uni-dimensional by \scite{OS},
\scite{STAB.5}, \scite{UNI}, and
there are $\lambda > \aleph_0,M^* \in K^c_1$,
ch$(M^*) = \lambda > \aleph_0,M^*$ is not $(D,\lambda)$-homogeneous.
By \scite{UNI.1}, $M^*$ is not $(D,\aleph_1)$-homogeneous.
\enddemo
\bn
Then there exists $B \subseteq M^*, \; |B| \le \aleph_0, \; p \in
\bold S^1_D(B),p$ is omitted in $M^*$, and in fact by \scite{MOM.4C},
$p^{[2 \varepsilon]}$ is omitted in $M^*$ for some $\varepsilon 
> 0$.  Therefore (by \scite{4n.9.1}) $p$ has no pseudo
$(<\varepsilon)$-support in $M^*$, see Definition \scite{4n.9}.

Let $\mu > \lambda$ be large enough.
\nl

We choose $a_\alpha$ by induction on $\alpha < \mu$ such that
\mr
\item "{$\circledast$}"  $A_\alpha = M^* \cup \{\bar a_\beta:\beta <
\alpha\}$ really $(<\varepsilon)$-omits $p(\bar x)$.
\endroster
\bn
\ub{Case (a)}:  If $A_\alpha =: M^* \cup \{\bar a_\beta:\beta <
\alpha\}$ is not in $K_1$.  Then first choose $\varphi(\bar x,\bar y) \in
\Delta,\bar a \subseteq A_\alpha$ such that $\varphi(\bar x,\bar a)$ 
witnesses $A_\alpha  \notin
K_1$, and second choose $\bar a_\alpha$ realize some strictly isolated $q \in
\bold S^{\ell g(\bar x)}_D(A_\alpha)$ which contains $\varphi(\bar
x,\bar a)$.

By \scite{8n.2} below, i.e., next claim, this is possible and
$\circledast$ is preserved.
\bn
\ub{Case (b)}: Not (a), then $\bar a_\alpha \notin A_\alpha$ and
$\circledast$ holds, using \scite{8n.3} below.
\nl
Having carried out the construction, let $M 
= \text{\rm mcl}(A_\mu)$.  On the one hand, $M$ belongs to $K^c_1$ 
as $A_\mu \in K_1$ (in fact, $A_\mu \prec {\frak C}$) by case (a) of
the construction.  On the other hand, note that $M$ is not 
$(D,\aleph_1)$-homogeneous, and moreover,
cannot be $\varepsilon$-isomorphic to a $(D,\aleph_1)$-homogeneous
model, as $p(x)$ is $(< \varepsilon)$-omitted by it ($p$ is really $(<
\varepsilon)$-omitted by $\circledast$, and recall
\scite{4n.9.1}(2)).  So we can construct such a model of arbitrarily
large density character, and by a L\"owenheim-Skolem argument,of
\ub{any} density character large enough.

But for each $\mu > \aleph_0$ there is $M \in K^c_1$, 
Ch$(M) = \mu$, $M$ is $(D,\aleph_1)$-homogeneous (by \scite{EM.1}(2)), a
contradiction to wu-categoricity.
\bn
In order to complete the proof of the main theorem, we only need to
show that the construction above is possible, which is done in the
following two claims.
\proclaim{\stag{8n.2} Claim}   [$({\frak C},\bold d)$ is
$0^+-\aleph_0$-stable, compact]

Let $p \in \bold S^1_D(B),B$ \ub{countable}, $A \supseteq B,A$ is not a
pseudo $(< \varepsilon)$-support for $p$, see Definition \scite{4n.9}. 
Let $\varphi(\bar x,\bar a)$ be a consistent formula over $A$.  \ub{Then} there
exists $\bar b \in {\frak C}$ such that ${\frak C} \models \varphi
(\bar b,\bar a)$ and $A \cup \bar b$ is not a pseudo 
$(< \varepsilon)$-support for $p$.  In fact, it is enough to choose
$\bar b$ such that {\rm tp}$(\bar b,A)$ is strictly
isolated and ${\frak C} \models \varphi(\bar b,\bar a)$.
\endproclaim
\bigskip

\demo{Proof}  By \scite{4n.12}(2) for some $\bar b \in \varphi({\frak
C},\bar a)$, tp$(\bar b,A)$ is strictly isolated.
So assume toward contradiction
\mr
\item "{$(*)_1$}"  $A \cup \bar b$ is a pseudo $(<
\varepsilon)$-support for $p(\bar x)$.
\ermn
Hence (by Definition \scite{4n.9}) there are
$\zeta(1),\zeta(2),\vartheta_1(\bar x,\bar b,\bar c)$
\mr
\item "{$(*)_2$}"  $\vartheta^{<\zeta(1)>}_1(\bar x,\bar b,\bar c_1)
\models p^{[\varepsilon-\zeta(2)]}(\bar x)$ and $\bar c_1 \subseteq
A$ and $\vartheta_1(\bar x,\bar b,\bar c_1)$ is consistent.
\ermn
As tp$(\bar b,A)$ is strictly isolated (see Definition \scite{4n.13}),
there are $\zeta(3) > 0$ and
$\psi(\bar y,\bar c_2)$ such that
\mr
\item "{$(*)_3$}"  $(i) \quad \psi(\bar y,\bar c_2) \in \text{\rm
tp}(\bar b,A)$ so $\bar c_2 \subseteq A$
\sn
\item "{${{}}$}"  $(ii) \quad \psi(\bar y,\bar c_2)$ pseudo
$(\zeta(1),\zeta(3))$-isolates tp$(\bar b,A)$,
i.e. $\psi^{<\zeta(3)>}(\bar y,c_2] \models \text{ tp}(\bar b,A)^{[\zeta(1)]}$.
\ermn
Let
\mr
\item "{$(*)_4$}"  $\vartheta_2(\bar x,\bar c_1,\bar c_2) = (\exists \bar
y)[\psi(\bar y,\bar c_2) \wedge \vartheta_1(\bar x,\bar y,c_1)]$.
\ermn
Clearly
\mr
\item "{$(*)_5$}"  $\vartheta_2(\bar x,\bar c_1,\bar c_2)$ is consistent.
\ermn
Choose $\zeta(4)$ such that
\mr
\item "{$(*)_6$}"  $0 < \zeta(4) < \zeta(3)$ and $\zeta(4) < \zeta(1)$.
\ermn
We shall now show that $\vartheta_2(x,\bar c_1,\bar c_2)$  
pseudo $(\varepsilon -\zeta(2),\zeta(4))$-isolates $p$ 
(and is over $A$), i.e.
\mr
\item "{$\boxtimes$}"  $\vartheta^{<\zeta(4)>}_2(\bar x,\bar c_1,\bar
c_2) \models p^{[\varepsilon-\zeta(2)]}(\bar x)$.
\ermn
This will give a contradiction to the assumption that $A$ is not a
pseudo $(< \varepsilon)$-support for $p$.
So assume
\mr
\item "{$(*)_7$}"  ${\frak C} \models \vartheta^{<\zeta(4)>}_2[\bar
a,\bar c_1,\bar c_2]$.
\ermn
By the definition of $\vartheta^{<\zeta(4)>}_2$, there are $\bar a',\bar
c'_1,\bar c'_2$ such that
\mr
\item "{$(*)_8$}"  $(i) \quad {\frak C} \models \vartheta_2[\bar a',
\bar c'_1,\bar c'_2]$
\sn
\item "{${{}}$}"  $(ii) \quad \bold d(\bar a',\bar a) \le \zeta(4) \le
\zeta(1)$
\sn
\item "{${{}}$}"  $(iii) \quad \bold d(\bar c'_1,\bar c_1) \le
\zeta(4) \le \zeta(1)$
\sn
\item "{${{}}$}"  $(iv) \quad \bold d(\bar c'_2,\bar c_2) \le
\zeta(4) \le \zeta(3)$.
\ermn
By the choice of $\vartheta_2$, i.e. $(*)_4$, for some $\bar b'$
\mr
\item "{$(*)_9$}"  $(i) \quad {\frak C} \models \psi[\bar b',\bar c'_2]$
\sn
\item "{${{}}$}"  $(ii) \quad {\frak C} \models
\vartheta_1[\bar a',\bar b',\bar c'_1]$.
\ermn
By the choice of $\psi(\bar y,\bar c_2)$, i.e., $(*)_3(ii)$ (note that 
as $0 < \zeta(4) \le \zeta(3)$, we have $\psi^{<\zeta(3)>}(\bar
b',\bar c_2))$ there is $\bar b'' \in {\frak C}$ such that
\mr
\item "{$(*)_{10}$}"  $(i) \quad \bar b''$ realizes tp$(\bar b,A)$
\sn
\item "{${{}}$}"  $(ii) \quad \bold d(\bar b'',\bar b') \le \zeta(1)$.
\ermn
So $\bold d(\bar a,\bar a') \le \zeta(4) \le \zeta(1),\bold d(\bar
b',b'') \le \zeta(1)$ and $\bold d(\bar c_1,\bar c'_1)
\le \zeta(4) \le \zeta(1)$, therefore by $(*)_9(ii)$ we get
\mr
\item "{$(*)_{11}$}"  ${\frak C} \models \vartheta^{<\zeta(1)>}_1
[\bar a,\bar b'',\bar c_1]$.
\ermn
So by $(*)_2$, replacing $\bar b$ with $\bar b''$ which has the same
type over $A$
\mr
\item "{$(*)_{12}$}"  $\bar a$ realizes $p^{[\varepsilon-\zeta(2)]}$
\ermn
so we have finished proving $\boxtimes$, hence getting the desired
contradiction.   \hfill$\square_{\scite{8n.2}}$
\enddemo
\bigskip

\proclaim{\stag{8n.3} Claim}  [$({\frak C},\bold d)$ is
$0^+-\aleph_0$-stable,uni-dimensional]
\nl
1) Assume 
\mr
\item "{$(a)$}"  $M = \text{\rm mcl}(M) \subseteq {\frak C}$
\sn
\item "{$(b)$}"  {\rm Ch}$(M) > \aleph_0$
\sn
\item "{$(c)$}"  $B \subseteq M$ is countable
\sn
\item "{$(d)$}"  $p \in \bold S^m_D(B)$
\sn
\item "{$(e)$}"  $M$ really $(< \varepsilon)$-omits $p(\bar x)$, see
Definition \scite{4n.9}.
\ermn
\ub{Then} for some $\bar b \in {\frak C} \backslash M$, also $M \cup
\{\bar b\}$ pseudo $(< \varepsilon)$-omits $p$.
\nl
2) Assume clauses (a)-(e) and $\bar b \in {\frak C},\bar b
 \notin M$, then $M \cup \bar b$ pseudo $(< \varepsilon)$-omits
$p(x)$ when: for every $\bar c \in {}^{\omega >}M$ and $\zeta >0$
there is $\bar b'$ realizing {\rm tp}$(\bar b,B \cup \bar c)$ such that
$\bold d(\bar b',{}^{\ell g(\bar b)}M) < \zeta$.
\nl
3) In clause (2) it is enough to assume that for every $\bar c \in
{}^{\omega >}M$ and $\zeta > 0$ there exist $\bar b',\bar b^-$ such
that $\bold d(\bar b,\bar b^-) < \zeta,\bold d(\bar b',M) < \zeta$
 and $\bar b'$ realizes tp$(\bar b^-,B \cup \bar c)$.
\endproclaim
\bigskip

\remark{Remark}  Assuming just $|\tau_{\frak C}| + \aleph_0 + |B| < \text{\rm
Ch}(M)$ suffices.  Assuming Ch$(M) > |B| + 2^{\aleph_0} + 
|\tau_{\frak C}|$, we can waive the assumption on $0^+-\aleph_0$-stability.
\endremark
\bigskip

\demo{Proof}  1) Let $\Phi^*$ be as in \scite{EM.1}(3), so
$\tau(\Phi^*)$ of cardinality $\aleph_1$.  Recall
\mr
\item "{$(*)$}"  for every uncountable (large enough if we do not
assume $\tau$ countable) linear order $I,M^*_I = \text{\rm mcl}
(\text{\rm EM}_\tau(I,\Phi^*))$ is $(D,|I|)$-homogeneous.
\ermn
[It is $(D,\aleph_1)$-homogeneous and use uni-dimensionality.]  
Choose $I=\lambda$ for $\lambda$ regular,
 $\lambda > \text{ Ch}(M)$, and let $M^* = M^*_\lambda$.

So \wilog \, $M \subseteq M^*$.  As Ch$(M) > \aleph_0$, we can by
\scite{C.2.B} below 
find $\varepsilon > 0$ and $a_\alpha \in M$ for $\alpha < \omega_1$
such that $\alpha < \beta < \omega_1 \Rightarrow \bold
d(a_\alpha,a_\beta) > \varepsilon$.  Let $\{a_{\omega_1+n}:n <
\omega\} \subseteq M$ include $B$.  Now for each $\alpha < \omega_1 +
\omega$ and $n <\omega$ we can find $b_{\alpha,n} \in \text{\rm
EM}_\tau(I,\Phi^*)$ such that $\bold d(a_\alpha,b_{\alpha,n}) <
1/(n+1)$.  Let $b_{\alpha,n} =
\sigma_{\alpha,n}(a_{t_{\alpha,n,0}},\dotsc,a_{t_{\alpha,n,k(\alpha,n)-1}})$.

Let (for $\alpha < \omega_1$)

$$
S_\alpha = \{t_{\beta,n,\ell}:(\beta < \alpha) \vee \beta \in
[\omega_1,\omega_1 +\omega) \text{ and }
n < \omega,\ell < k(\beta,n)\} \subseteq \lambda.
$$
\mn
Let (for $\alpha < \omega_1,n < \omega,\ell < k(\alpha,n))$

$$
\gamma_{\alpha,n,\ell} = \text{\rm Min}\{\gamma \in S_\alpha \cup
\{\lambda\}:t_{\alpha,n,\ell} \le \gamma\}.
$$
\enddemo
\bigskip

\demo{\stag{8n.3.1} Subclaim}  Under these assumptions, there exists
$C$,  a club of $\omega_1$, such that
\mr
\item "{$(*)$}"  if $\delta \in C$ and $m < \omega$ then the following
set is stationary
$$
\align
W_{\delta,m} = \{\alpha \in C:&\text{for every } n \le
m,\sigma_{\alpha,n} = \sigma_{\delta,n}, \text{ hence } k(\alpha,n) =
k(\delta,n) \text{ and} \\
  &(\gamma_{\alpha,n,\ell} = \gamma_{\delta,n,\ell}) \wedge
  (\gamma_{\alpha,n,\ell} \in S_\alpha \equiv \gamma_{\delta,n,\ell}
  \in S_\delta) \text{ for } \ell < k(\delta,n)\}.
\endalign
$$
\endroster
\enddemo
\bigskip

\demo{Proof}   By a standard coding argument, there exist
functions $f_n:\omega_1 \rightarrow \omega_1$ (for $n < \omega$) such
that for $\alpha < \omega_1,f_n(\alpha)$ ``encodes" the finite sequence

$$
\langle \sigma_{\alpha,m}:m \le n\rangle \frown
\langle(\beta^*_{\alpha,m,\ell},m^*_{\alpha,m,\ell},\ell^*_{\alpha,m,\ell}):m
\le n,\ell < k(\alpha,m)\rangle
$$
\mn
where $\beta^*_{\alpha,m,\ell} < \omega_1 +
\omega,m^*_{\alpha,m,\ell}$ and $\ell^*_{\alpha,m,\ell}$ are natural
numbers, satisfying:
\mr
\item "{$\otimes$}"
$(\beta^*_{\alpha,m,\ell},m^*_{\alpha,m,\ell},\ell^*_{\alpha,m,\ell})$
is the minimal (lexicographically) triple $(\beta^*,m^*,\ell^*)$ such
that $t_{\beta^*,m^*,\ell^*} = \gamma_{\alpha,m,\ell}$.
\ermn
In fact, there exists such coding $f_n:\omega_1 \rightarrow \omega_1$ 
such that on a club $C'_n,f_n$ is regressive.  Now by F\"odor's lemma,
the following set contains a club (as its $\omega_1$ complement cannot
contain a stationary set):

$$
C''_n = \{\delta \in C'_n:\{\alpha \in C'_n:f_n(\alpha) =
f_n(\delta)\} \text{ is stationary}\}.
$$
\mn
We call this club $C_n$ and let $C = \dbca_{n < \omega} C_n$, obviously
$C$ is as required.  \hfill$\square_{\scite{8n.3.1}}$
\enddemo
\bigskip

\demo{\stag{8n.3A} Subclaim}  Under these assumptions, 
there exist $\delta(*) \in C$ 
such that $(\forall m < \omega)(\exists \delta \in
\delta(*) \cap C)(\delta(*) \in W_{\delta,m})$.
\enddemo
\bigskip

\demo{Proof}  Note that $\delta(*) \in W_{\delta,m}
\Leftrightarrow \delta \in W_{\delta(*),m} \Leftrightarrow
W_{\delta,m} = W_{\delta(*),m}$, so all we are looking for is
$\delta(*)$ satisfying $\delta(*) > \text{ min }W_{\delta(*),m}$ for
all $m$; now if for all $\delta < \omega_1 \exists m < \omega$ such
that $\delta = \text{ min } W_{\delta,m}$ we set an easy
contradiction (e.g. by F\"odor's lemma, although it is an overkill here).
\hfill$\square_{\scite{8n.3A}}$ 

By a similar argument we can find $\delta(*) \in C$ such that there
exists a sequence $\langle \delta_n:n < \omega\rangle,\delta_n \in
C,\delta_n < \delta_{n+1} < \delta(*)$ for all $n$ and
$W_{\delta(*),n} = W_{\delta_n,n}$ for all $n$, so in particular $m
\le n \Rightarrow W_{\delta_m,m} = W_{\delta_n,m} = W_{\delta(*),m}$.  
\nl
We obtain (by the choice of $\delta_n,C,W_{\delta_n,m}$, etc.):
\mr
\item "{$\circledast_1$}"  $\langle b_{\delta_n,m}:n \ge m\rangle$ is
an indiscernible sequence over $[B_m = B \cup\{b_{\alpha,i}:\alpha <
\delta_m,i < \omega\}]$ for each $m$.
\ermn
In fact, we can say more:
\mr
\item "{$\circledast_2$}"  for each $\bar c \in \text{ EM}(I,\Phi^*)$
finite and for each $m < \omega$, there exists $n^* < \omega$ such that 
$\langle b_{\delta_n,m},n \ge n^*\rangle$ is indiscernible 
over $B_m \cup \bar c$.
\ermn
We would like now to continue the indiscernible sequences above in a
proper extension of $M^+$.

Let $J = I \times \Bbb Q$ and let us identify $I$ with $I \times
\{0\}$, so we think of $I$ as a subset of $J$.  Now 
for each $m < \omega,\ell < \kappa(
\delta_m,m)$, look at the sequence $\langle t_{\delta_n,m,\ell}:n \ge
m\rangle$.  It is either constant or strictly increasing (recall
$\circledast_1$), in the first case define $t^*_{m,\ell} = t_{\delta_m,m,\ell}
\in I \subseteq J$, otherwise choose $t^*_{m,\ell} \in J$ such that
$t_{\delta_n,m,\ell} < t^*_{m,\ell} < \text{\rm
sup}\{t_{\delta_n,m,\ell}:n \ge m\}$.

Note:
\mr
\item "{$(**)$}"  for each $m < \omega$, the (quantifier free) type
of the sequence $\langle t^*_{n,\ell}:n < m,\ell < k(\delta(*),n
\rangle$ in the language $\{<\}$ (order) is the same as the type of
$\langle t_{\delta(*),n,\ell}:n < m,\ell < k(\delta(*),n)\rangle$.
\ermn
Let $M^+ = \text{\rm mcl}(\text{\rm EM}_{\tau({\frak C})}(J,\Phi^*))
\prec^1_\Delta {\frak C},M^+$ extends $M^*$.  Let $b^*_n =
\sigma_{\delta(*),n}(t^*_{n,0},\dotsc,t^*_{n,k(\delta(*),n)-1})$ so $b^*_n
\in M^+,\langle b^*_n:n < \omega \rangle$ is a Cauchy sequence 
(as $\langle b_{\delta^*,n}:n < \omega \rangle$ is a Cauchy sequence,
use the indiscernibility of $J$) with limit $b^* \in M^+$.

Note (by $\circledast_2$ and the choice of $b^*_n$)
\mr
\item "{$\circledast_3$}"  for each $\bar c \in \text{ EM}(I,\Phi^*)$
finite and for each $m < \omega$ there exists $n^* < \omega$ such that
$\langle b_{\delta_n,m}:n \ge n^*\rangle \frown b^*_m$ is
indiscernible over $B \cup \bar c$.  In fact, by \scite{MON.2s.1}, the
same is true for each $\bar c \in M^*$ finite.
\ermn
We also observe
\mr
\item "{$\circledast_4$}"   $b^* \in M^+ \backslash M^*$.
\ermn
Why?  Recall that $\langle a_\alpha:\alpha <
\omega_1\rangle$ form an $\varepsilon$-net, so for $m$ big enough,
$b_{\delta_{n_1},m},b_{\delta_{n_2},m}$ can not be too close for $n_1
\ne n_2$.  Combining this with $\circledast_3$ we see that $b^*_m$ can
not be $\frac{\varepsilon}{3}$ close to any $\bar c \in M^*$.
\mr
\item "{$\circledast_5$}"   $b^*$ is as required in (3).
\ermn
Why?  Given $\bar c \in M$ finite and $\zeta >0$ we choose $m$ and $n$
big enough such that $b^*_m$ is close enough to $b^*,b_{\delta_n,m}$
 is close enough to $a_{\delta_n} \in M_1$, and tp$(b_{\delta_1,m},B
\cup \bar c) = \text{ tp}(b^*_m,B \cup \bar c)$ [possible by $\circledast_3$].
\nl

2) Assume this fails, so there are $\zeta > 0,\bar c \subseteq M$ and a
formula $\vartheta(\bar x,\bar b,\bar c)$ 
such that $\vartheta^{<\zeta>]}(x,\bar
b,\bar c) \models p^{[\varepsilon-\zeta]}(\bar x)$.

Choose $\zeta_1,\zeta_2 > 0$ such that $\zeta_1 < \zeta_2 < \zeta$ and
$\zeta_1 < \zeta - \zeta_2$.  By the assumption there are $\bar
b',\bar b''$ such that
\mr
\item "{$(*)$}"  $(a) \quad \bar b''$ realizes tp$(\bar b,B \cup \bar c)$
\sn
\item "{${{}}$}"  $(b) \quad \bar b' \subseteq M$ of length $\ell
g(\bar b)$
\sn
\item "{${{}}$}"  $(c) \quad \bold d(\bar b',\bar b'') < \zeta_1 <
\zeta - \zeta_2$.
\ermn
By clause (a) of $(*)$ as $\theta^{<\zeta>}(\bar x,\bar b,\bar c)
\models p^{[\varepsilon-\zeta]}(\bar x)$ also $\theta^{<\zeta>}(\bar
x,\bar b'',\bar c) \models p^{[\varepsilon-\zeta]}(\bar x)$.  Hence
(by $*(c)$ and \scite{4n.10}(3)) $\theta^{<\zeta_2>}(\bar x,\bar
b',\bar c) \models p^{[\varepsilon-\zeta]}(\bar x)$, a contradiction.
\nl
3) Similar proof to (2), first using \scite{4n.10}(3) to show that
$\theta^{<\zeta'>}(\bar x,\bar b^-,\bar c) \models 
p^{[\varepsilon-\zeta]}(\bar x)$ for some $\zeta'$.
\hfill$\square_{\scite{8n.3}}$ 
\enddemo
\bigskip

\demo{\stag{C.2.B} Observation}   Let 
$(X,\bold d)$ be a non-separable metric space.   Then there exist 

$$
\langle a_i:i < \omega_1\rangle \subseteq X \text{ and } \varepsilon^*
> 0 \text{ such that } \bold d(a_i,a_j) \ge \varepsilon^* \, \forall i,j <
\omega_1.
$$
\enddemo
\bigskip

\demo{Proof}
Choose by induction on $i < \omega_1,a_i$ such that 
$a_i \notin \{a_j:j < i\}$.   Choose 
$0 < \varepsilon_i \le d(a_i,\{a_j:j < i\}),\varepsilon_i \in \Bbb Q$.
Without loss of generality  $\varepsilon_i = \varepsilon^*$ for all
$i$, and we are done.  \hfill$\square_{\scite{CATEG}}$
\enddemo
\bigskip

\nocite{ignore-this-bibtex-warning} 
\newpage
    
REFERENCES.  
\bibliographystyle{lit-plain}
\bibliography{lista,listb,listx,listf,liste}

\enddocument